\documentclass[10.5pt,a4paper]{article}
\usepackage{amsfonts,amsmath,amsthm,amsxtra,amscd,amssymb}
\usepackage{amsmath}
\usepackage{amsthm}
\usepackage{amsxtra}
\usepackage{amsbsy}
\usepackage{bm}
\usepackage{newlfont}
\usepackage{pictexwd,dcpic}
\usepackage[dvips]{color}
\usepackage{bm}
\usepackage{newlfont}
\usepackage[dvips]{color}
\theoremstyle{plain}
\newtheorem{thm}{Theorem}[section]
\newtheorem{prop}[thm]{Proposition}
\newtheorem{lem}[thm]{Lemma}
\newtheorem{cor}[thm]{Corollary}

\theoremstyle{definition}
\newtheorem{exmp}{Example}[section]
\newtheorem{defn}{Definition}[section]
\newtheorem{rem}{Remark}[section]

\def\P #1{\partial_{#1}}
\def\ch #1{\text{\rm Char}\, #1}
\def\Q #1{\partial^2_{#1}}
\def\pl{{\partial}}

\def\G{\Gamma}
\def\a{\alpha}
\def\b{\beta}

\def\a{\alpha}
\def\b{\beta}

\def\G{\Gamma}

\def\L{\Lambda}
\def\O{\Omega}

\def\o{\omega}
\def\O{\Omega}

\def\s{\sigma}
\def\S{\Sigma}
\def\t{\tau}

\begin{document}
\title{Geometry of Darboux-Manakov-Zakharov Systems and its Application}
\author{Peter J. Vassiliou, \\ Faculty of Information Sciences and Engineering,\\ University of Canberra, A.C.T. AUSTRALIA, 2601\\ {\tt peter.vassiliou@canberra.edu.au}}                     
\maketitle
\begin{abstract}
The intrinsic geometric properties of generalized Darboux-Manakov-Zakharov systems of semilinear partial differential
equations
\begin{equation}\label{GDMZabstract}
\frac{\partial^2 u}{\partial x_i\partial x_j}=f_{ij}\Big(x_k,u,\frac{\partial u}{\partial x_l}\Big),\ \ \ 
1\leq i<j\leq n,\ \ k,l\in\{1,\ldots,n\}
\end{equation} 
for a real-valued function $u(x_1,\ldots,x_n)$ are studied with particular reference to the linear systems in this equation class.
 
System (\ref{GDMZabstract}) 
will not generally be involutive in the sense of Cartan: its coefficients will be constrained by complicated nonlinear integrability conditions. We derive geometric tools for explicitly constructing involutive systems of the form (\ref{GDMZabstract}), essentially solving the integrability conditions. Specializing to the linear case provides us with a novel way of viewing and solving the multi-dimensional $n$-wave resonant interaction system and its modified version as well as constructing new examples of semi-Hamiltonian systems of hydrodynamic type. The general theory is illustrated by a study of these applications. 
%
\vskip 10 pt
\noindent{\it Keywords:} Three-wave resonant interaction system, Darboux-Manakov-Zakharov, semi-Hamiltonian hydrodynamic type, local differential geometry
\vskip 5 pt
\noindent{\it 2010 Mathematics Subject Classification:} 35N05, 53B05, 58A17, 58A30, 58D19, 58J45
\end{abstract}
\section{Introduction}

This paper describes a geometrical study of the system of nonlinear partial differential equations
\begin{equation}\label{nWaveSystem}
\frac{\partial A_{j\,k}}{\partial x_i}=A_{j\,i}A_{i\,k},\ \ \ \ \text{for all}\ \ (i,\,j,\,k)\in \text{perm}_3(n) 
\end{equation}
where $\text{perm}_3(n)$ denotes the set of all permutations of 3 letters out of $n$. Note that the convention of summation over repeated indices is not in force in equation (\ref{nWaveSystem}) and will not be used in this paper. Any summation will make explicit use of the summation symbol.  

As we will explain in more detail in the next section, system (\ref{nWaveSystem}) plays a important role in aspects of submanifold geometry and integrable systems theory [\ref{Darboux10}, \ref{ZakharovManakov85}, \ref{Zakharov98}, \ref{Tsa}, \ref{JoshiKitaevTreharne08}] as well as arising as model equations in a variety of physical phenomena including nonlinear optics and plasma physics [\ref{ZakharovManakov73}, \ref{ZakharovManakov85}]. It also plays a  role in mathematical physics, for instance, 
[\ref{MangazeevSergeev01}, \ref{Sergeev07}].

A particularly important special case of (\ref{nWaveSystem}) is when $n=3$, known as the $2+1$-{\it dimensional three-wave resonant interaction system}
({\sf 3WRI}), which serves as a model equation for nonlinear wave interaction. 

The main objectives of this paper are firstly, to introduce a new geometric setting for system (\ref{nWaveSystem}), namely, {\it $n$-hyperbolic manifolds}, and then show how this leads to new solutions and new perspectives on this integrable system. Secondly, to describe some of the consequences arising from these results for the above mentioned areas of application especially to {\it systems of hydrodynamic type}, which have a rich mathematical structure. The central object that draws all these many issues together is a class of overdetermined systems of linear hyperbolic partial differential equations, the {\it Darboux-Manakov-Zakharov systems} which we now discuss.

\section{Darboux-Manakov-Zakharov Systems}
\label{intro}
To set the results to be presented in context, we provide background on the crucial tool that is used in this paper, namely, Darboux-Manakov-Zakharov systems and explain the links that they have with various nonlinear phenomena. 

The explicit modern study of overdetermined linear systems of partial differential equations of the form
\begin{equation}\label{KTsystem1}
\frac{\pl^2 u}{\pl x_i\pl x_j}-\G_{j\,i}(x)\frac{\pl u}{\pl x_i}-
\G_{ij}(x)\frac{\pl u}{\pl x_j}+C_{ij}(x)u=0,\   
\ 1\leq i< j\leq n,  
\end{equation}
for a real-valued function $u(x_1,\ldots,x_n)$, where\footnote{Warning: We take this opportunity to remind the reader that the convention of summation over repeated indices {\it will not} be used in this paper and is not in force in equation (\ref{KTsystem1}).} coefficients 
$\G_{ij},\G_{j\,i}, C_{ij}=C_{j\,i}$ are smooth functions of the $x_i$    
began in the work of Manakov \& Zakharov [\ref{ZakharovManakov85}, \ref{Zakharov98}] and Kamran \& Tenenblat [\ref{KT96}, \ref{KT98}, \ref{KT2000}]. The former studies were motivated by nonlinear optics and founded on the extensive classical works of 19th century geometers on the problem of triply orthogonal coordinate systems, particularly those of Darboux [\ref{Darboux1896}, \ref{Darboux10}]. The studies of Kamran and Tenenblat arose in their search for a suitable venue for generalising the classical Laplace Transformation to higher dimensions following work of Chern's [\ref{Chern44}] on the theory of Cartan submanifolds. The work of all these authors turns out to be extremely interesting and is connected with a variety of phenomena in submanifold geometry and integrable systems theory. 

However a challenge is that for $n>2$, system (\ref{KTsystem1}) is not involutive in general and the applications mentioned above very much rely on its involutivity. This state of affairs presents both a problem and an opportunity since the complicated integrability conditions that arise have been the object of study in mathematical physics for about three decades, at least in the case $n=3$, and much longer in differential geometry [\ref{Darboux1896}]. In the 3-dimensional case, they are equivalent to the well known 2+1-dimensional 3-wave resonant interaction ({\sf 3WRI}) equations and form part of the Lam\'e equations for triply orthogonal coordinate systems [\ref{ZakharovManakov85}, \ref{Zakharov98}, \ref{Darboux10}].  

Until now essentially two approaches have been taken, in the recent literature, to the construction of solutions of the 2+1-dimensional {\sf 3WRI} system. Firstly, the  Darboux-Manakov-Zakharov  linear problem [\ref{ZakharovManakov85}] or else the Ablowitz-Haberman linear problem [\ref{AblowitzHaberman75}] have been used to obtain some solutions via the inverse scattering transform (IST) [\ref{Cornille79}, \ref{Kaup80}, \ref{Kaup81}, \ref{Zakharov98}]. In relation to this we mention also the works [\ref{Kaup81a}, \ref{LeviPilloniSantini81}] which construct B\"acklund transformations for (among other equations) the 2+1-dimensional {\sf 3WRI} system. Important solutions of the 2+1-dimensional {\sf 3WRI} system obtained in this way include the lump solutions [\ref{Kaup81}]. However, it should be mentioned that by comparison to other integrable systems, the IST solution of the {\sf 3WRI} system is very complicated and not so staightforward to apply.

Secondly, Lie symmetry methods have been extensively applied to obtain explicit solutions; references here include [\ref{Kitaev90}, \ref{LeoMartinaSoliani86}, \ref{MartinaWinternitz89}]. In the last reference a detailed study of the Lie point symmetry group of the 2+1-dimensional {\sf 3WRI} system is shown to be infinite-dimensional with Lie algebra of Kac-Moody-Virosoro type. Numerous explicit group-invariant solutions are constructed via symmetry reduction using the one- and two-dimensional Lie subgroups of the full symmetry group of the system. In particular, as first pointed out by Kitaev [\ref{Kitaev90}], there are solutions expressible in terms of the Painlev\'e transcendents. This fact has been used in more recent work [\ref{JoshiKitaevTreharne08}] to study the linearization of Painlev\'e equations.

Of course, another source of solutions of the {\sf 3WRI} system and extensive information on linear system (\ref{KTsystem1}) comes to us from the nineteenth century geometers. Many of these results are summarised in Darboux's famous treatise [\ref{Darboux10}].   Additionally, formulas relevant to the modern treatment of the Lam\'e equations via the inverse scattering transform are discussed in detail in the classical work [\ref{Darboux1896}, Chapter 12], without however, explicitly anticipating IST itself.  A linear PDE system of the form (\ref{KTsystem1}) {\it which is involutive} in the sense of Cartan will be called a {\it linear Darboux-Manakov-Zakharov system} which we abreviate to {\it linear {\sf DMZ} system} or simply {\it {\sf DMZ} system} whenever there is little danger of confusion.

This paper takes a radically different approach to the problem of constructing solutions of the 2+1-dimensional {\sf 3WRI} equations by providing a geometric characterisation and construction of {\it generalised Darboux-Manakov-Zakharov systems} 
\begin{equation}\label{GDMZ}
\frac{\partial^2 u}{\partial x_i\partial x_j}=f_{ij}\Big(x_k,u,\frac{\partial u}{\partial x_l}\Big),\ \ \ 
1\leq i<j\leq n,\ \ k,l\in\{1,\ldots,n\},
\end{equation}
which we sometimes denote {\it {\sf GDMZ} systems}, and then showing how to specialize to linear {\sf DMZ} systems. In consequence new tools are created for the study of linear {\sf DMZ} systems with remarkably many applications arising from the central role they play in the theory of integrable equations, submanifold geometry and mathematical physics. For instance, we construct, in case $n=3$, linear {\sf DMZ} systems which therefore provide new explicit solutions for the 2+1-dimensional {\sf 3WRI} equations. Note that the restriction to the case $n=3$ is not essential so that all the results can be easily extended to  $n>3$.  It is noteworthy that the problem of constructing linear {\sf DMZ} systems was first proposed by Darboux [\ref{Darboux1896}, Ch. 12]. There Darboux derives numerous results subject to the proviso of explicit knowledge of a linear {\sf DMZ} system; he states\footnote{pp 278-279}
\begin{center}
{\it ``Supposons que l'on ait obtenu par un moyen quelconque un syst\'eme\\  d'\'equations aux  d\'eriv\'ees partielles} { [\,(\,\ref{KTsystem1}\,)\,]} {\it pour lequel les conditions d'int\'egrabilit\'e soient satisfaites et dont on ait d\'etermin\'e l'int\'egrale $u$."}  
\end{center}
Darboux's motivation for this problem arose from the fact that linear {\sf DMZ} systems are deeply linked to the problem of constructing $n$-orthogonal coordinate systems in an $n$-dimensional flat Riemannian manifold, say $\mathbb{E}^n$, $n$-dimensional Euclidean space. To see this, let $x^1,\ldots,x^n$ be coordinates on $\mathbb{E}^n$. The problem of $n$-orthogonal coordinate systems is to find functions $u^i(x^1,\ldots,x^n)$, $i=1,\ldots,n$ such that the map $\boldsymbol{u}:x\mapsto u(x)$ is a local diffeomorphism and satisfies the PDE system
$$
\sum_{k=1}^n\frac{\partial u^i}{\partial x^k}\frac{\partial u^j}{\partial x^k}=0,\ \ \ \ \forall\ \ \ \ i,j\in\{1,2,\ldots,n\},\ \ \ i\neq j.
$$ 
The problem for $n=2$ is easy but much more difficult for $n\geq 3$ and it was intensively studied by numerous important 19th century geometers, leading to the so called  {\it Lam\'e equations}. We will briefly outline the construction of these equations as they will play a role in subsequent sections of this paper. Let us denote the inverse of $\boldsymbol{u}$ by
$
\boldsymbol{x}(u^1,\ldots,u^n)
$
and define
$$
h_i^2=\sum_{k=1}^n\left(\frac{\partial x^i}{\partial u^k}\right)^2.
$$
The metric in $\mathbb{E}^n$ expressed in the $(u^i)$-coordinate system has the form
$$
\sum_{i=1}^n\,h_i^2\,(du^i)^2.
$$
The Christoffel symbols for this metric are
$$
\begin{aligned}
&\Gamma^i_{ml}=0,\ \ \ \ \ i\neq l\neq m,\cr
&\Gamma^i_{il}=\frac{1}{h_i}\frac{\partial h_i}{\partial u^l},\cr
&\Gamma^i_{ll}=-\frac{h_l}{h_i^2}\frac{\partial h_l}{\partial u^i},\ \ \ \ \ \ i\neq l.
\end{aligned}
$$
The vanishing of the Riemann curvature tensor leads to the PDE system
\begin{equation}\label{LameEquations}
\begin{aligned}
&\frac{\partial^2 h_i}{\partial u^l\partial u^m}-
\frac{1}{h_l}\frac{\partial h_l}{\partial u^m}\frac{\partial h_i}{\partial u^l}-
\frac{1}{h_m}\frac{\partial h_m}{\partial u^m}\frac{\partial h_i}{\partial u^l}=0,\cr
&\frac{\partial}{\partial u^l}\left(\frac{1}{h_l}\frac{\partial h_i}{\partial u^l}\right)
+\frac{\partial}{\partial u^i}\left(\frac{1}{h_i}\frac{\partial h_l}{\partial u^i}\right)
+\sum_{k\neq i,l}^n\,\frac{1}{h_k^2}\frac{\partial h_i}{\partial u^k}\frac{\partial h_l}{\partial u^k}=0
\end{aligned}
\end{equation}
One can show that system (\ref{LameEquations}) is involutive in the sense of the Cartan-K\"ahler theorem [\ref{BC3G},\ref{IveyLandsberg03}] and the solution depends upon $n(n-1)/2$ functions, each of 2 variables. 

We will now make the connection between the $2+1$-dimensional {\sf 3WRI} equations, the Lam\'e equations and {\sf DMZ} system (\ref{KTsystem1}) more precise by formulating the following result that is essentially contained in [\ref{Darboux1896}, Chapter 12].

\begin{thm}[Darboux]\label{geometricBackgroundThm}
Suppose the linear PDE system
\begin{equation}\label{3WRI_LP}
\frac{\partial^2u}{\partial x_i\partial x_j}-
\G_{j\,i}(x)\frac{\partial u}{\partial x_i}-\G_{ij}(x)\frac{\partial u}{\partial x_j}=0,\ \ \ \ 1\leq i<j\leq 3
\end{equation}
for functions $u:\mathbb{R}^3\to\mathbb{R}$, where $\G_{ij},\G_{j\,i}$ are functions of $x=(x_1,x_2,x_3)$ is involutive. Then
\begin{enumerate}
\item[{\rm 1.}] There are real-valued functions $h_i(x)$, $1\leq i\leq 3$, such that
$$
\G_{ij}=\frac{1}{h_j}\frac{\partial h_j}{\partial x_i}
$$
\item[{\rm 2.}] If $u$ is a solution of (\ref{3WRI_LP}) then the functions
\begin{equation}\label{3WRI_def}
\psi_i=\frac{1}{h_i}\frac{\partial u}{\partial x_i},\ \ \ \    
\ A_{ij}=\frac{1}{h_i}\frac{\partial h_j}{\partial x_i}
\end{equation}
satisfy the 2+1-dimensional 3-wave resonant interaction system and its linear problem 
\begin{equation}\label{3WRI}
\frac{\partial A_{j\,k}}{\partial x_i}=A_{j\,i}A_{i\,k}\ \ \ \ \ 
\frac{\partial \psi_j}{\partial x_i}=A_{j\,i}\psi_i,
\end{equation}
respectively
\item[{\rm 3.}] The functions $h_1,h_2,h_3$ satisfy the nonlinear system
\begin{equation}\label{halfLame}
\frac{\partial^2 h_i}{\partial x_j\partial x_k}-
\frac{1}{h_j}\frac{\partial h_j}{\partial x_k}\frac{\partial h_i}{\partial x_j}-
\frac{1}{h_k}\frac{\partial h_k}{\partial x_k}\frac{\partial h_i}{\partial x_j}=0,
\ \ \ (i,j,k)\in\ \text{\rm perm}\,(1,2,3).
\end{equation}
\end{enumerate}
\end{thm}

\proof{  The integrability conditions for linear system (\ref{3WRI_LP}) are
\begin{equation}\label{KT_0Compatibility}
\frac{\partial\G_{ij}}{\partial x_k}-\G_{ik}\G_{kj}-\G_{ki}\G_{ij}+\G_{ij}\G_{kj}= 
\frac{\partial \G_{ij}}{\partial x_k}-\frac{\partial \G_{kj}}{\partial x_i}=0,\ \ (i,j,k)\in\ \text{\rm perm}\,(1,2,3)
\end{equation}
From the second set of equations in (\ref{KT_0Compatibility}) we deduce that there are real valued functions functions $h_j$ on $\mathbb{R}^3$ such that
$$
\G_{ij}=\frac{1}{h_j}\frac{\partial h_j}{\partial x_i}
$$
Next, let $u$ be any solution of (\ref{3WRI_LP}) and define functions $\psi_i$ and $A_{ij}$ as in (\ref{3WRI_def}). Then we have
$$
\begin{aligned}
\frac{\pl \psi_j}{\pl x_i}-A_{ji}\psi_i=&
\frac{\pl}{\pl x_i}\left(\frac{1}{h_j}\frac{\pl u}{\pl x_j}\right) -\frac{1}{h_j}\frac{\pl h_i}{\pl x_j}\frac{1}{h_i}
\frac{\pl u}{\pl x_i}\cr
=&\frac{1}{h_j}\left( \frac{\pl^2 u}{\pl x_i\pl x_j}
-\frac{1}{h_i}\frac{\pl h_i}{\pl x_j}\frac{\pl u}{\pl x_i}
-\frac{1}{h_j}\frac{\pl h_j }{\pl x_i}\frac{\pl u}{\pl x_j} \right)\cr
=&\frac{1}{h_j}\left(\frac{\pl^2 u}{\pl x_i\pl x_j}
-\G_{j\,i}\frac{\pl u}{\pl x_i}-\G_{ij}\frac{\pl u}{\pl x_j}\right)=0.
\end{aligned} 
$$ 
It is then easily verified that the equations $(\ref{3WRI})_2$ have the 2+1-dimensional 3-wave resonant interaction equations 
$(\ref{3WRI})_1$ as integrability conditions. This proves  assertion 2 of the theorem. Finally, substituting the form of $A_{ij}$ from (\ref{3WRI_def}) into $(\ref{3WRI})_1$ gives (\ref{halfLame}).}\hfill \qed  

\vskip 5 pt
Notice that equations (\ref{halfLame}) form part of the Lam\'e equations for triply orthogonal coordinate
systems, namely, $(\ref{LameEquations})_2$. When the $h_i$ satisfy only (\ref{halfLame}) rather than system (\ref{LameEquations}), they determine a triply orthogonal coordinate system in a Riemannian manifold $M$ which is not flat. However, it can be shown that the metric of $M$
$$
g=\sum_{i=1}^n h_i^2\,dx_i^2
$$
in the coordinates $x_i$ that arises from any solution of (\ref{halfLame}) has diagonal Ricci tensor. Moreover, the coefficients $\G_{ij},\, \G_{ji}$ in {\sf DMZ} system (\ref{3WRI_LP}) agree with the Levi-Civita connection coefficients $\G^i_{ij},\, \G^j_{ij}$, for $g$, respectively. 
 
Zakharov [\ref{Zakharov98}] calls such Riemannian manifolds {\it spaces of diagonal curvature} and compares them with systems of hydrodynamic type which are integrable in a sense to be described but which do not possess a Hamiltonian structure. These are the so called {\it semi-Hamiltonian systems} about which we shall say more below.  

Note also that equations (\ref{3WRI}) 
contain the 2+1-dimensional 3-wave resonant interaction equations, $(\ref{3WRI})_1$, together with it's standard linear problem $(\ref{3WRI})_2$. In other words, system $(\ref{3WRI})_1$ expresses the integrability conditions of  $(\ref{3WRI})_2$, in the linear case (\ref{3WRI_LP}).

From this theorem we see that the linear system (\ref{3WRI_LP}) is the key that unifies the 2+1-dimensional {\sf 3WRI} system and certain triply orthogonal coordinate systems {\it provided (\ref{3WRI_LP}) is involutive}. An important goal of this paper is to derive feasible means of constructing linear {\sf DMZ} systems (\ref{3WRI_LP}) without being required to directly solve the nonlinear integrability conditions (\ref{KT_0Compatibility}) for the coefficients $\G_{ij},\G_{j\,i}$ in (\ref{3WRI_LP}). In fact, we examine the {\it geometric properties} of generalised {\sf DMZ} systems (\ref{GDMZabstract}) in order to characterise and construct them, thereby constructing new solutions for the 2+1-dimensional {\sf 3WRI} system.

As a concrete example of a {\sf DMZ} system, consider the triply orthogonal system defined by {\it oblate spheroidal coordinates}
$$
u=\cosh\,x\,\,\cos\,y\,\,\cos\,z,\ \ \ v=\cosh\,x\,\,\cos\,y\,\,\sin\,z,\ \ \
w=\sinh\,x\,\,\sin\,y.
$$
The corresponding flat metric is
$$
g=(\cosh^2x-\cos^2y)\,(dx^2+dy^2)+\cosh^2x\ \cos^2y\ dz^2
$$
and the relevant linear system corresponding to (\ref{3WRI_LP}) is
$$
\frac{\pl^2 u }{\pl x_i\pl x_j}-\G^i_{ij}\frac{\pl u}{\pl x_i}
-\G^j_{ij}\frac{\pl u}{\pl x_j}=0.
$$
Computing Christoffel coefficients delivers the linear system 
\begin{equation}\label{oblateSpheroidal}
\begin{aligned}
&\frac{\pl^2 u }{\pl x\pl y}-\left(\frac{\cos\,y\ \sin\,y}{\cosh^2x-\cos^2y}\right)\frac{\pl u}{\pl x}
-\left(\frac{\cosh\,x\ \sinh\,x}{\cosh^2x-\cos^2y}\right)
\frac{\pl u}{\pl y}=0,\cr
&\frac{\pl^2 u }{\pl x\pl z}-\tanh\,x\frac{\pl u}{\pl z}=0,
\ \ \frac{\pl^2 u }{\pl y\pl z}+\tan\,y\frac{\pl u}{\pl z}=0.
\end{aligned}
\end{equation}
System (\ref{oblateSpheroidal}) is involutive in the sense of the Cartan-K\"ahler theorem, as can be verified; it is therefore a {\sf DMZ} system. One way to see this is to restrict the contact system on jet space $J^2(\mathbb{R}^3,\mathbb{R})$
to the 10-dimensional submanifold $\mathcal{R}\subset J^2(\mathbb{R}^3,\mathbb{R})$ defined by (\ref{oblateSpheroidal}). This defines a rank 4 linear Pfaffian system $\theta^i$, $0\leq i\leq 3$ on $\mathcal{R}$ with independence forms $\o^1=dx,\o^2=dy,\o^3=dz$. These 7 differential 1-forms can be extended to a coframe on $\mathcal{R}$ by $\pi_1,\pi_2,\pi_3$ satisfying the structure equations

\begin{equation}\label{mainStructure}
\begin{aligned}
&d\theta^0\equiv 0,\cr
&d\theta^1\equiv \pi_1\wedge\o^1,\cr
&d\theta^2\equiv \pi_2\wedge\o^2,\cr
&d\theta^3\equiv \pi_3\wedge\o^3,
\end{aligned}\mod\{\theta^0,\theta^1,\theta^2,\theta^3\}.
\end{equation} 
 
This verifies that there is no essential torsion, that is, no integrability conditions, and that the tableau [\ref{BC3G},\,\ref{IveyLandsberg03}] satisfies Cartan's involutivity test.  Accordingly, by Theorem \ref{geometricBackgroundThm}, the functions $A_{ij}$ defined by 
(\ref{3WRI_def}) satisfy the 2+1-dimensional {\sf 3WRI} system. However, there are many more {\sf DMZ} systems (\ref{3WRI_LP}) than those provided by triply orthogonal coordinate systems and in this paper we will show how to construct many such systems. In consequence of Theorem 1.1, we construct orthogonal coordinate systems in spaces of diagonal curvature.

An integrable system of partial differential equations closely related to the {\sf 3WRI} system is the 2+1-dimensional {\it modified three-wave resonant interaction system}, which we denote {\sf m3WRI}, discovered by V. Mangazeev and S. Sergeev [\ref{MangazeevSergeev01}, \ref{Sergeev07}]. Like the {\sf 3WRI} system, the {\sf m3WRI} system has a Hamiltonian structure and plays a role in discrete differential geometry and quantum field theory [\ref{Sergeev07}]; it has the form
\begin{equation}\label{m3WRI}
\frac{\pl \G_{jk}}{\pl x_i}=(\G_{ij}-\G_{ik})(\G_{jk}-\G_{ji}),\ \ \ (i,j,k)\in\text{perm}\{1, 2, 3\}.
\end{equation}    
As shown by Mangazeev and Sergeev, it arises as the integrability condition of a {\sf DMZ} system of the form
\begin{equation}\label{MangazeevSergeev_LP}
\left(\P {x_i}\P {x_j}-\G_{ij}\P {x_j}-\G_{ji}\P {x_i}+
\G_{ij}\G_{ji}\right)\phi=0,\ \ \ 1\leq i<j\leq 3.
\end{equation}
We demonstrate in this paper how our general construction gives rise to new solutions of the {\sf m3WRI} system by constructing {\sf DMZ} systems of the form (\ref{MangazeevSergeev_LP}).

The problem of orthogonal coordinates and the 2+1-dimensional {\sf 3WRI} system are, in turn, linked to certain nonlinear partial differential equations of mathematical physics that have attracted a good deal of interest in recent years. To explain this let $\boldsymbol{u}=(u^1,u^2,\ldots,u^n)$ be a vector valued function of space variable $x$ and time $t$ and consider PDE of the form
\begin{equation}\label{quasilinear}
\frac{\partial u^i}{\partial t}-\sum_{k=1}^n\,v^i_k(\boldsymbol{u})\frac{\partial u^k}{\partial x}=0.
\end{equation} 
Dubrovin and Novikov [\ref{DN}] showed that 
if there is a flat Levi-Civita connection $\nabla$ and function $h$ on phase space such that
$$
v^i_j=\nabla^i\nabla_j\,h,
$$
then (\ref{quasilinear}) has a Hamiltonian structure. Subsequently, S. Tsarev showed that among non-Hamiltonian systems (\ref{quasilinear}) there is a significant subset of them which possess remarkable integrability properties.
To describe these we restrict ourselves to systems (\ref{quasilinear}) that are expressible in distinct Riemann invariants; these are the {\it strongly hyperbolic} systems. For strongly hyperbolic systems there is a change of dependent variables such that 
$$
(v^j_k)=\text{diag}(w^1(\bar{u})\ \ w^2(\bar{u})\ \ \ldots\ \ w^n(\bar{u}))
$$
and for all $i,j$, $w^i(\bar{u})\neq w^j(\bar{u})$.
Attention focusses on those strongly hyperbolic systems that satisfy the additional constraint
\begin{equation}\label{diagonalCondition}
\partial_k\left(\frac{\partial_j w^i}{w^j-w^i}\right)=\partial_j\left(\frac{\partial_k w^i}{w^k-w^i}\right).
\end{equation}

Even though strongly hyperbolic systems (\ref{quasilinear}) whose ``diagonal velocities" satisfy (\ref{diagonalCondition}) are not, in general Hamiltonian, they nevertheless have significant integrability properties; for this reason they are called {\it semi-Hamiltonian}. For such systems Tsarev [\ref{Tsa}] has established a number of notable properties (see section 5 of this paper). The main point to note at this stage is that any {\sf DMZ} system (\ref{3WRI_LP}) gives rise to a semi-Hamiltonian PDE system in (1+1)-dimensions (\ref{quasilinear}).  We demonstrate in section 5 how the geometric properties of {\sf DMZ} systems derived in sections 2 and 3, lead to new examples of such systems. 

Finally we mention a link between the {\sf DMZ} systems (\ref{KTsystem1}) and a class of submanifolds whose study was initiated by Cartan [\ref{Cartan19}] and then subsequently developed by Chern [\ref{Chern44}] and much more recently by Kamran and Tenenblat in their generalisation of Laplace transformations to dimensions higher the two [\ref{KT96}, \ref{KT98}, \ref{KT2000}]. Below we schematically summarise the central role that {\sf DMZ} systems  (\ref{KTsystem1}) play in the above mentioned fields of submanifold geometry and integrable systems.

\vskip 5 pt

\begin{equation}
\begindc{\commdiag}[40]
\obj(0,0){DMZ}
\obj(1,0){E}
\obj(-1,0){W}
\obj(0,-1){S}
\obj(0,1){N}
\mor{DMZ}{N}{}
\mor{DMZ}{S}{}
\mor{DMZ}{E}{}
\mor{DMZ}{W}{}
\enddc
\end{equation}

where
\begin{enumerate}
\item \underline{DMZ} denotes the Darboux-Manakov-Zakharov linear systems (\ref{KTsystem1}) in the case $n=3$; [\ref{ZakharovManakov85}, \ref{Zakharov98}].
\item \underline{N} denotes Lam\'e Equations for Triply Orthogonal Coordinate Systems; [\ref{Zakharov98}]
\item \underline{S} denotes Semi-Hamiltonian Systems of Hydrodynamic Type; [\ref{Tsa}], [\ref{Se}]
\item \underline{E} denotes the 2+1-dimensional Three-Wave Resonant Interaction Systems; [\ref{Kaup80},\ref{Sergeev07}]
\item \underline{W} denotes Cartan Sub-manifolds and Higher-Dimensional Laplace Transformations; [\ref{KT96}]
\end{enumerate}

\vskip 5 pt

We will have opportunity, in sections 4 and 5 of this paper, to discuss and illustrate the links \underline{S} and \underline{E} in more detail, using the theory developed in sections 3 and 4. For the present, it suffices to emphasise the importance of finding a means of characterising and constructing {\sf GDMZ} systems (\ref{GDMZ}) since it will enable us to significantly affect a number of problems in geometry, integrable systems and mathematical physics. We turn to this problem next.

\section{Intrinsic geometry of {\sf GDMZ} partial differential equations}
\label{sec:1}

In this section we derive a method for geometrically characterising involutive overdetermined systems of partial differential equations in 1 dependent variable $u$ and $n$ independent variables $x_1,\ldots,x_n$ of the form
\begin{equation}\label{mainEquations}
\frac{\partial^2 u}{\partial x_i\partial x_j}=f_{ij}\Big(x_k,u,\frac{\partial u}{\partial x_l}\Big),\ \ \ 1\leq i<j\leq n,\ \ k,l\in\{1,\ldots,n\}.
\end{equation}
In particular, we will be interested in {\it explicitly} constructing {\it linear} systems in the class (\ref{mainEquations}). 
\begin{defn}
In case system (\ref{mainEquations}) is involutive, then we shall refer to it as a {\it generalised Darboux-Manakov-Zakharov} ({\sf GDMZ}) PDE system. If the $f_{ij}$ in (\ref{mainEquations}) are such that it is a {\rm linear} involutive PDE system then we shall refer to it as a {\it linear Darboux-Manakov-Zakharov} system or simply as a {\it Darboux-Manakov-Zakharov} ({\sf DMZ}) system.  
\end{defn}
For $n=2$ system (\ref{mainEquations}) comprises just one equation and is automatically involutive. For $n>2$ the system is overdetermined; it's involutive if and only if integrability conditions are satisfied. The contribution of this section is to derive a geometric characterisation of systems (\ref{mainEquations}) for which the integrability conditions are satisfied ``automatically". That is, a geometric characterisation of generalised Darboux-Manakov-Zakharov systems. A significant by-product of this result in the linear case, is the construction of explicit solutions of the 3-wave resonant interaction equations and the modified 3-wave resonant interaction equations. This will be discussed in section 4. Another is the construction of new semi-Hamiltonian systems of hydrodynamic type, explained in section 5.

The abovementioned geometric characterisation of {\sf GDMZ} systems relies on the existence of certain structured manifolds and special coordinate systems on them. For any $n\geq 2$, let $M_{3n+1}$ be a smooth manifold of dimension $3n+1$ equiped with a rank $2n$ distribution $H$ which decomposes as a direct sum of $n$ rank 2 distributions $H_i$
\begin{equation}\label{nHyp}
H=\bigoplus_{i=1}^nH_i
\end{equation}
satisfying the structure equations
\begin{equation}\label{genHypStr}
[H_i,H_j]\equiv 0,\ [H_i,H_i]\equiv Z_i\mod H,\ \forall\,i\ \text{and}\ j\neq i
\end{equation}
where $Z_1\wedge\cdots\wedge Z_n\neq 0$. Establishing that such distributions exist in plentiful supply, can be explicitly constructed and have interesting properties and applications are aims of this paper.

\begin{defn}
Let $n\geq 2$ be an integer. A pair $(M_{3n+1},H)$ consisting of a manifold $M$ of dimension $3n+1$ and a distribution $H$ of the form (\ref{nHyp}) satisfying structure equations  (\ref{genHypStr}) and derived type 
$$
\mathfrak{d}(H)=[[2n,0],[3n,n],[3n+1,3n+1]]
$$
will be called an $n$-{\it hyperbolic manifold}, while $H$ will be called an $n$-{\it hyperbolic structure} or $n$-{\it hyperbolic distribution}.\hfill $\clubsuit$
\end{defn}

Note that the $j^{\,\text{th}}$, ordered, 2-component element of the list of lists $\mathfrak{d}(H)$ records the dimension $\dim\,H^{(j)}$ of the $j^{\,\text{th}}$ derived distribution $H^{(j)}$ of $H$ and the dimension $\dim\,\ch\,H^{(j)}$ of its Cauchy distribution 
$$
[\dim\,H^{(j)},\dim\,\ch\,H^{(j)}].
$$
Recall that for any smooth distribution $\mathcal{V}\subseteq TM$ on a smooth manifold $M$, the {\it Cauchy} or {\it characteristic} distribution, $\ch\mathcal{V}$ is defined by
$$
\ch\mathcal{V}=\{X\in\mathcal{V}~|~[X,\mathcal{V}]\subseteq\mathcal{V}\}.
$$
The derived distribution $\mathcal{V}^{(1)}$ of $\mathcal{V}$ is defined by
$$
\mathcal{V}^{(1)}=\pl\mathcal{V}=\mathcal{V}+[\mathcal{V},\mathcal{V}]
$$
and the higher derived distributions are defined by iteration: $\mathcal{V}^{(j)}=\pl^j\mathcal{V}$. It is easy to prove, in the case $\mathcal{V}$ and $\ch\mathcal{V}$ have constant rank on $M$ that $\ch\mathcal{V}$ is Frobenius integrable. See [\ref{Vassiliou06a},\, \ref{Vassiliou06b}], for further details. Throughout this paper, constant rank assumptions are made for all vector field distributions and codistributions.

\begin{defn}
Let $n\geq 2$ be an integer. We say that an $n$-hyperbolic manifold $(M_{3n+1},H)$ is {\it Daboux integrable} if there is a coordinate system
\begin{equation}\label{adaptedCoords}
\text{\sf C}: x_1,x_2,\ldots,x_n,y_0,y_1,\ldots,y_n,u_1,\ldots,u_n
\end{equation}
on $M$,  {\it adapted} to its $n$-hyperbolic structure $H$ in the sense that
\begin{enumerate}
\item[1).]
The Cauchy distribution ${\ch H^{(1)}}$ of the derived distribution $H^{(1)}$ satisfies 
$$
\dim \ch H^{(1)}=n,\ \ \ch H^{(1)}\subset H
$$ 
and 
$x_1,x_2,\ldots,x_n,y_0,y_1,\ldots,y_n$ are the invariants of $\ch H^{(1)}$, while $u_1,\ldots,u_n$ extend these to a coordinate system on $M$ from which it follows that
$$
\ch H^{(1)}=\{\P {u_1},\ldots,\P {u_n}\}
$$
\item[2).]  For each $i$, there is a basis for $H$ of the form
\begin{equation}\label{adaptedBasis}
X^1_i=\P {x_i}+P_i,\ \ X^2_i=\P {u_i},\ \ 1\leq i\leq n
\end{equation}
where
$$
P_i=\sum_{a=0}^n\rho^a_i(x,y,u)\P {y_a}\mod\{\P {u_1},\ldots,\P {u_n}\}.
$$
and $\{X^1_i,\P {u_i}\}$ is a basis for $H_i$, for some real-valued functions $\rho^a_i$.
\end{enumerate}
The vector fields (\ref{adaptedBasis}) of a Darboux integrable $n$-hyperbolic manifold are said to form an {\it adapted basis} for $H$.\hfill $\clubsuit$
\end{defn}
That Darboux integrable $n$-hyperbolic manifolds exist and can be constructed will be established in section 3. 
Our aim now is to show that a Darboux integrable $n$-hyperbolic manifold determines an immersion $\Phi$ of the $M$ into $J^2(\mathbb{R}^n,\mathbb{R})$ whose image is a  {\sf GDMZ} system, (\ref{mainEquations}). 

From now on we drop the superscript 1 on $X^1_i$ and rename this vector field $X_i$. 

\begin{lem}\label{AdaptedCoords1}
Let $(M_{3n+1},H)$ be an Darboux integrable $n$-hyperbolic manifold with adapted coordinates (\ref{adaptedCoords}) and adapted basis $\{X_i,\P {u_i}\}$ for $H_i$. Then for all $i\neq j$ with $i,j\in \{1,2,\ldots,n\}$ we have
$$
[X_i,X_j]\equiv 0,\ [X_i,\P {u_j}]\equiv 0\mod\{\P {u_1}, \P {u_2},\ldots,\P {u_n}\}.
$$
\end{lem}
\vskip 5 pt
\proof{ By the hyperbolic structure equations (\ref{genHypStr}) we have 
$$
[X_i,X_j]=\sum_{k=1}^nA^k_{ij}X_k+\sum_{k=1}^nB^k_{ij}\P {u_k}
$$ 
for some functions $A,B$. However, the assumed local form (\ref{adaptedBasis}) of the adapted basis for $H$ shows that $A^k_{ij}=0$ for all $i, j, k$. A similar argument establishes the second congruence.}\hfill\qed

\vskip 5 pt

\begin{lem}\label{AdaptedCoords2}
Let $(M,H)$ be a Darboux integrable $n$-hyperbolic manifold with adapted coordinates (\ref{adaptedCoords}) and adapted basis $\{X_i,\P {u_i}\}$ for $H_i$. Then 
$$
\mathcal{A}=\Big\{\P {u_1},\P {u_2},\ldots,\P {u_n},[\P {u_1},X_1],[\P {u_2},X_2],\ldots,[\P {u_n},X_n]\Big\}\subset H^{(1)}
$$
is a rank $2n$, integrable distrubution.
\end{lem}
\vskip 5 pt
\proof{
Since $\dim H^{(1)}=3n$ and the $X_j$ pairwise commute modulo $\{\P {u_1}\ldots,\P {u_n}\}$ by the previous Lemma, it follows that $\dim\mathcal{A}=2n$. By the Jacobi identity and again by the previous Lemma, we deduce 
for all $i\neq j$ that $[\P {u_j},[\P {u_i},X_i]]\equiv 0\mod \{\P {u_1}\ldots,\P {u_n}\}$. Similarly, we deduce that $[X_j,[\P {u_i},X_i]]\equiv 0\mod \{\P {u_1}\ldots,\P {u_n}\}$.  Next set $\xi_i=[\P {u_i},X_i]$, apply the Jacobi identity to get $[\xi_i,[\P {u_j},X_j]]+[X_j,[\xi_i,\P {u_j}]]+[\P {u_j},[X_j,\xi_i]]=0$. From this and using our previous results we deduce that 
$$
[[\P {u_i},X_i],[\P {u_j},X_j]]\equiv 0\mod\mathcal{A}.
$$
Hence $\mathcal{A}$ is rank $2n$ and integrable.}\hfill\qed

\begin{cor}\label{AdaptedCoords3}
Let $(M,H)$ be an Darboux integrable $n$-hyperbolic manifold. The invariants of the distribution $\mathcal{A}\subset H^{(1)}$ are spanned by $x_1,\ldots,x_n,p$ for some function $p$ on $M$.  
\end{cor}

\begin{thm}\label{KTdistributions}
Let $(M,H)$ be a Darboux integrable $n$-hyperbolic manifold. Let functions $x_1,\ldots,x_n,p$ on $M$ span the invariants of $\mathcal{A}$ and for $1\leq i<j\leq n$, define functions $p_i,p_{ii},f_{ij}$ by
$$
p_i=X_ip,\ p_{ii}=X_ip_i,f_{ij}=X_jp_i.
$$
Then the map $\varphi$ given by
$$
\varphi(x,y,u)=(x_1,\ldots,x_n,p,p_1,\ldots,p_n,p_{11},\ldots,p_{nn})
$$
is a local diffeomorphism and defines an immersion
$$
\Phi : M\hookrightarrow J^2(\mathbb{R}^3,\mathbb{R})
$$
whose image is a {\sf GDMZ} system
$$
\frac{\partial^2 u}{\partial x_i\partial x_j}=f_{ij}\Big(x_k,u,\frac{\partial u}{\partial x_l}\Big),\ \ \ 
1\leq i<j\leq n,\ \ k,l\in\{1,\ldots,n\}.
$$
\end{thm}

\vskip 5 pt

\proof First we show that the functions
$$
x_1,\ldots,x_n,\, p,\, p_1,\ldots,p_n,\, p_{11},\ldots,p_{nn}
$$
as defined in the statement of the theorem form a local coordinate system on $M$. Begin by extending the independent invariants $x_1,\ldots,x_n,p$ of $\mathcal{A}$ by $y_1,\ldots,y_n,\allowbreak u_1,\ldots,u_n$, where the $y_l$ extend $x_1,x_2\ldots,x_n,p$ to be a complete set of invariants of $\ch H^{(1)}$. In these coordinates, the local form of each $X_i$ is
$$
X_i=\P {x_i}+\rho^0_i\P p+\sum_{l=1}^n\rho^l_i\P {y_l}\mod \P {u_j}
$$
and we therefore have
$$
[\P {u_i},X_i]=\frac{\partial \rho^0_i}{\partial u_i}\P p+
\sum_{l=1}^n\frac{\partial \rho^l_i}{\partial u_i}\P {y_l}\mod \{\P {u_j}\},\ \ 1\leq i\leq n.
$$
But the left-hand-side is an element of $\mathcal{A}$ and hence
$$
0=[\P {u_i},X_i]p=\frac{\partial \rho^0_i}{\partial u_i}
$$
and thus
$$
[\P {u_i},X_i]=
\sum_{l=1}^n\frac{\partial \rho^l_i}{\partial u_i}\P {y_l}\mod \{\P {u_j}\},\ \ 1\leq i\leq n.
$$
It follows that the matrix
$$
J=\left(\frac{\partial\rho^j_i}{\partial u_i}\right)_{i,j=1}^n
$$
is nonsingular. We deduce that
$$
H^{(1)}=\Big\{\P {x_1}+\rho^0_1\P p,\ldots,\P {x_n}+\rho^0_n\P p,\P {y_1}, \ldots, \P {y_n},\P {u_1}, \ldots, \P {u_n}\Big\}
$$
and that the quotient of $H^{(1)}$ by $\ch H^{(1)}$
$$
K:=H^{(1)}/\ch H^{(1)}=\Big\{\P {x_1}+\rho^0_1\P p,\ldots,\P {x_n}+\rho^0_n\P p,\P {y_1}, \ldots, \P {y_n}\Big\}
$$
is a rank $2n$ distribution defined on the $2n+1$-dimensional quotient manifold\newline 
$M/\ch H^{(1)}$, satisfying $\dim K^{(1)}=2n+1$
and $\ch{K}=\{0\}$. It follows from Pfaff-Darboux theorem that new coordinates
$$
p_l=\rho^0_l
$$
can be taken in place of the $y_l$, expressing $K$ in the canonical form
$$
K=\Big\{\P {x_1}+p_1\P p,\ldots,\P {x_n}+p_n\P p,\P {p_1}, \ldots, \P {p_n}\Big\}.
$$
This means there are functions $\widetilde{\rho}^{\,l}_i$ of the independent coordinates on $M$ given by $x,p,p_l,u_j$ such that
$$
X_i=\P {x_i}+p_i\P p+\sum_{l=1}^n\widetilde{\rho}^{\,l}_i\P {p_l}\mod\{\P {u_j}\}.
$$
The functions $p_{ij}$ are defined by
\begin{equation}\label{derivatives2}
p_{ij}=X_jp_i=\widetilde{\rho}^{\,j}_i.
\end{equation}
and since the $X_i$ pairwise commute modulo $\ch\,H^{(1)}$, we have
\begin{equation}\label{lower_symmetric}
p_{ij}=\widetilde{\rho}^{\,j}_i=\widetilde{\rho}^{\,i}_j=p_{ji}.
\end{equation}
To see this we calculate
$$
p_{ij}-p_{j\,i}=X_ip_j-X_jp_i=[X_i,X_j]p=0
$$
since $[X_i,X_j]\in\{\P {u_1},\, \P {u_2}\ldots,\P {u_n}\}$, for all $i,j$, by Lemma 2.1.

Furthermore, we have
\begin{equation}\label{1stDerivatives}
\begin{aligned}
\frac{\partial p_k}{\partial u_l}&=\frac{\partial}{\partial u_l}\left(X_kp\right)\cr
&=\frac{\partial}{\partial u_l}\left(X_kp\right)-X_k\left(\frac{\partial p}{\partial u_l}\right),\ \text{by}\ \text{Corollary}\ \ref{AdaptedCoords3}\cr
&=[\P {u_l},X_k]p=0
\end{aligned}
\end{equation}
because $[\P {u_l},X_k]\in\mathcal{A}$ for all $k,l$ by Lemmas \ref{AdaptedCoords2} and \ref{AdaptedCoords3}.  This means that
$$
[\P {u_l},X_k]=\frac{\pl p_k}{\pl u_l}\P p+\sum_{i=1}^n\frac{\pl\widetilde{\rho}^{\,i}_k}{\pl u_l}\P {p_i}=
\sum_{i=1}^n\frac{\pl\widetilde{\rho}^{\,i}_k}{\pl u_l}\P {p_i}.
$$
But for $l\neq k$, we have $[\P {u_l}, X_k]\in\{\P {u_1},\ldots,\P {u_n}\}$ and hence
\begin{equation}\label{nondiagonal}
\frac{\partial \widetilde{\rho}^{\,i}_k}{\partial u_l}=0,\ \ \ \ \forall\ \ i\ \text{and}\ l\neq k.
\end{equation}
Equations (\ref{derivatives2}), (\ref{lower_symmetric}) and (\ref{nondiagonal}) imply that the matrix $J$ has the form
$$
J=\text{diag}\,\left(\begin{matrix} \frac{\partial\rho^1_1}{\partial u_1} & \frac{\partial\rho^2_2}{\partial u_2} &
\ldots & \frac{\partial\rho^n_n}{\partial u_n}\end{matrix}\right)
$$
and since $J$ is non-singular we deduce that
\begin{equation}\label{indepDerivatives2}
\prod_{i=1}^n\frac{\partial\rho^i_i}{\partial u_i}\neq 0
\end{equation}
The condition (\ref{indepDerivatives2}) and our argument till now proves the claim that the functions defined by the theorem statement are functionally independent.

Next, we prove that the local diffeomorphism 
$$
\varphi: (x,y,u)\mapsto (x_1,\ldots,x_n,p,p_1,\ldots,p_n,p_{11},\ldots,p_{nn})
$$
induces an immersion of $M$ into $J^2(\mathbb{R}^3,\mathbb{R})$ whose image is a {\sf GDMZ} system. We will prove that $\varphi$ identifies $H$ with the canonical distribution on 
$J^2(\mathbb{R}^3,\mathbb{R})$ spanned by
$$
\Big\{D_1^f,\ldots,D_n^f,\P {p_{11}},\ldots,\P {p_{nn}}\Big\}
$$  
where
$$
D_i^f=\P {x_i}+p_i\P p+\sum_{j=1}^np_{ij}\P {p_j},
$$
and where for all $i\neq j$
$$
p_{ij}=p_{j\,i}=f_{ij}.
$$
More precisely, we prove that there is a function $\boldsymbol{g} : M\to GL(2n)$ such that
\begin{equation}
\varphi_*\mathbb{X}=\boldsymbol{g}\;\mathfrak{D}
\end{equation}
where
\begin{equation}\label{equivalence}
\begin{aligned}
&\mathbb{X}=\left(\begin{matrix} X_1 & X_2 &\cdot&\cdot & X_n &\P {u_1} &\P {u_2}
&\cdot &\cdot &\P {u_n}\end{matrix}\right)^T,\cr
&\mathfrak{D}=\left(\begin{matrix} D^f_1 & D^f_2 &\cdot &\cdot & D^f_n &\P {p_{11}} &
\P {p_{22}} &\cdot &\cdot &\P {p_{nn}}\end{matrix}\right)^T.
\end{aligned}
\end{equation}
In fact we will show that $\boldsymbol{g}$ has the form
\renewcommand{\arraystretch}{1.3}
\begin{equation}\label{G_structure}
\boldsymbol{g}=\left(\begin{matrix} I_n & \mathbf{g}^1_n\cr
                     \mathbf{0}_n &\mathbf{g}^2_n\end{matrix}\right)
\end{equation}
where $I_n,\mathbf{0}_n$ are the $n\times n$ identity and zero matrices respectively; $\mathbf{g}^a_n$ are certain $n\times n$ matrices determined by the equivalence $\varphi$. Recalling that the components of local diffeomorphism $\varphi$
are $x_i,p,p_j,p_{kk}$, then from (\ref{equivalence}) and (\ref{G_structure}) we obtain the equations
\begin{equation}\label{equivalencePDE}
\begin{aligned}
&p_i=X_ip,\ \ \P {u_i}p=0,\ 1\leq i\leq n\cr
&p_{ii}=X_ip_i,\ \ f_{ij}=X_jp_i,\ \ \ 1\leq i<j\leq n,\cr 
&\P {u_k}p_l=0,\ \ \ \ k,l\in\{1,2,\ldots,n\}.
\end{aligned}
\end{equation} 
for the components of $\varphi$. The remaining equations define the entries of $\mathbf{g}^1_n$ and $\mathbf{g}^2_n$.
Equations $(\ref{equivalencePDE})_1$ and $(\ref{equivalencePDE})_2$ are satisfied by definition of $p,\,p_i$ and $p_{ii}$.  Equations $(\ref{equivalencePDE})_3$ are satisfied by virtue of (\ref{1stDerivatives}).  Hence, we have shown that $\varphi$ maps the distribution $H$ to the canonical distribution for the overdetermined system
$$
\frac{\partial^2 u}{\partial x_i\partial x_j}=f_{ij}\Big(x_k,u,\frac{\partial u}{\partial x_l},
\frac{\partial^2 u}{\partial x_l^2}\Big),\ \ \ 
1\leq i<j\leq n,\ \ k,l\in\{1,\ldots,n\}.
$$
Next we show that the $f_{ij}$ do not depend upon the second order partial derivatives
$$
\frac{\partial^2 u}{\partial x_l^2}.
$$
That is, we show that $f_{ij}$ don't, in fact, depend upon $p_{ll}$ for any $l$. For suppose for some $i,j$,\ $i\neq j$ and some $l$, that 
$$
\frac{\partial f_{ij}}{\partial p_{ll}}\neq 0.
$$  
We calculate
$$
\begin{aligned}
\frac{\partial f_{ij}}{\partial u_l}&=
\sum_{m=1}^n\left(\frac{\partial f_{ij}}{\partial x_m}\frac{\partial x_m}{\partial u_l}+
\frac{\partial f_{ij}}{\partial p}\frac{\partial p}{\partial u_l}+
\frac{\partial f_{ij}}{\partial p_m}\frac{\partial p_m}{\partial u_l}+
\frac{\partial f_{ij}}{\partial p_{mm}}\frac{\partial p_{mm}}{\partial u_l}\right)\cr
&=\frac{\partial f_{ij}}{\partial p_{ll}}\frac{\partial p_{ll}}{\partial u_l}\ \ \ \ \text{(no summation over $l$})\cr
&\neq 0\ \ \ \ \ \ \text{by equation}\ \ (\ref{indepDerivatives2})
\end{aligned}
$$
However, this contradicts equation (\ref{nondiagonal}). Hence we have shown that $\varphi$ maps $H$ to the canonical distribution for the PDE system  
\begin{equation}\label{MyEquations}
\frac{\partial^2 u}{\partial x_i\partial x_j}=f_{ij}\Big(x_k,u,\frac{\partial u}{\partial x_l}\Big),\ \ \ 
1\leq i<j\leq n,\ \ k,l\in\{1,\ldots,n\}.
\end{equation}
as we wanted. 

Finally, we prove that the $n$-hyperbolic structure $H$, and hence equation (\ref{MyEquations}), is involutive with respect to the independence form
$dx_1\wedge dx_2\wedge\cdots\wedge dx_n$. In particular, the integrability conditions for (\ref{MyEquations}) are satisfied.  Set $\Theta=H^\perp$. We will show that $\Theta$ is an involutive linear Pfaffian system with respect to independence form $\o^1\wedge \o^2\wedge\cdots\wedge\o^n$, where $\o^j=dx_j$. Recall the structure equations for $H$ have the form
\begin{equation}\label{StructureEqs1}
\begin{aligned}
&[X_i, \P {u_i}]\equiv Z_i\mod H,\ \ \  1\leq i\leq n,\ \ \ \ \text{and},\cr
&[X_i, X_j]\equiv 0,\ \ \ [X_i, \P {u_j}]\equiv 0\mod \ch H^{(1)},\ \ i,j\in\{1,2,\ldots,n\},\ i\neq j.
\end{aligned}
\end{equation}
where $Z_1\wedge Z_2\cdots\wedge Z_n\neq 0$. Consequently, we can write
$$
H^{(1)}=\{X_1,\ldots,X_n,\P {u_1},\ldots,\P {u_n}\}\oplus\{Z_1,\ldots,Z_n\}
$$
which we can extend by $U$ to make a frame on $M$. We let
\begin{equation}
\o^1,\o^2,\ldots,\o^n,\pi^1,\pi^2,\ldots,\pi^n,\theta^0,\theta^1,\ldots,\theta^n
\end{equation}
be the coframe dual to 
$$
X_1,\ldots,X_n,\P {u_1},\ldots,\P {u_n}, U, Z_1,\ldots,Z_n.
$$
In particular $\{\theta^0,\theta^1,\ldots,\theta^n\}=H^\perp=\Theta$. If $X,Y$ are vector fields on $M$ then use of the identity
$$
d\theta^p(X,Y)=X\left(\theta^p\rfloor Y\right)-X\left(\theta^p\rfloor Y\right)-\theta^p\,\rfloor\,[X,Y]
$$
and structure equations (\ref{StructureEqs1}) permits one to establish the structure equations
\begin{equation}\label{formsStructureEqs}
\begin{aligned}
&d\theta^0\equiv 0,\cr
&d\theta^1\equiv \pi^1\wedge\o^1,\cr
&d\theta^2\equiv \pi^2\wedge\o^2,\cr
&\ \ \ \ \ \ \vdots\cr
&d\theta^n\equiv \pi^n\wedge\o^n
\end{aligned}\mod\Theta
\end{equation}
Equations (\ref{formsStructureEqs}) prove that $\Theta$ has no integrability conditions (zero essential torsion) and an involutive tableau with one nonzero Cartan character, $s_1=n$. 
\hfill\qed

\begin{prop}[Linear {\sf DMZ} systems]
Let
\begin{equation}\label{linearDMZfull}
\frac{\pl^2 u}{\pl x_i\pl x_j }-\sum_{k=1}^n\G^k_{ij}(x)\frac{\pl u}{\pl x_k}+C_{ij}(x)u=0,\ \ \ 1\leq i<j\leq n,
\end{equation}
be a system of linear partial differential equations for one dependent variable $u$ in independent variables $x=(x_1,\ldots x_n)$, $n\geq 2$, where $\G^k_{ij},\ C_{ij}$ are real vauled smooth functions, symmetric in their lower indices. Then system (\ref{linearDMZfull}) is involutive if and only if 
$$
\G^k_{i j}=0,\ \ \ \ \text{for distinct}\ i, j, k, 
$$ 
and for each $(i,j,k)\in\text{\rm perm}_3(n)$,
\begin{equation}\label{DMZfullCompatibility}
\begin{aligned}
&\frac{\pl\G^j_{ij}}{\pl x_k}-\G^k_{ik}\G^j_{kj}-\G^i_{ki}\G^j_{ij}+\G^j_{ij}\G^j_{kj}+C_{ik}=0,\cr 
&\frac{\pl \G^j_{ij}}{\pl x_k}-\frac{\pl \G^j_{kj}}{\pl x_i}=0,\cr
&\frac{\pl C_{ij}}{\pl x_k}-\frac{\pl C_{ik}}{\pl x_j}+
C_{kj}(\G^k_{ik}-\G^j_{ij})+C_{ij}\G^i_{ik}-C_{ik}\G^i_{ij}=0.  
\end{aligned}
\end{equation}
If (\ref{linearDMZfull}) is involutive then its solutions depend upon $n$ functions, each of one variable.
\end{prop}
\vskip 3 pt
\proof This is a calculation using Cartan's theory of linear Pfaffian systems, [\ref{BC3G}, \ref{IveyLandsberg03}], applied to (\ref{linearDMZfull}). The left hand sides of (\ref{DMZfullCompatibility}) are components of the torsion tensor. The top reduced Cartan character is $s_1=n$. \hfill\qed

\vskip 5 pt

The main point here is that system (\ref{DMZfullCompatibility}) is precisely the {\sf nWRI} system in the case $C_{ik}=0$, as explained in Theorem \ref{geometricBackgroundThm}, and it is the modified {\sf nWRI} system in case $C_{ik}=\G^i_{ik}\G^k_{ik}$. Hence,  if we can explicitly construct Darboux integrable $n$-hyperbolic manifolds whose associated {\sf GDMZ} systems, as guaranteed by  Theorem \ref{KTdistributions}, are {\it linear} then we will have constructed 
\begin{enumerate}
\item[{\bf i)}] Solutions of the multi-dimensional {\sf nWRI} system if $C_{ik}=0$
\item[{\bf ii)}] Solutions of the modified {\sf nWRI} system if $C_{ik}=\G^i_{ik}\G^k_{ik}$
\item[{\bf iii)}] Semi-Hamiltonian systems of strongly hyperbolic conservation laws and their commuting flows if $C_{ik}=0$
\end{enumerate}

In sections 3 and 4, we will lay out theory that solves each of these construction problems. In sections 4 and 5, we illustrate these constructions with explicit applications, which show that the calculations require only the Frobenius theorem as well as quadrature, entailing integration which is very often easy to carry out. This procedure provides solutions for the full integrability conditions (\ref{DMZfullCompatibility}) in general and, in particular, leads to a very wide class of solutions of the multi-dimensional $n$-wave resonant interaction systems as well as the construction of new semi-Hamiltonian systems of hydrodynamic type and their commuting flows.

\section{Constructing a class of {\sf GDMZ} systems from   $n$-hyperbolic manifolds}

In the previous section we proved, in Theorem \ref{KTdistributions}, that to explicitly construct a {\sf GDMZ} system it is sufficient to have constructed an explicit Darboux integrable $n$-hyperbolic manifold. In this section a method will be described for constructing a large class of these. 

We hasten to point out that Darboux integrable $n$-hyperbolic manifolds are far more numerous than those which we construct in this section. However, in this paper we are primarily interested in {\it linear} {\sf DMZ} systems and for this purpose the construction described below is adequate.  Having said this, nonlinear {\sf GDMZ} systems are not excluded from our construction.

Consider cartesian products of jet spaces of the form 
\begin{equation}\label{jetProds}
\mathbb{J}=J^{k_1}(\mathbb{R},\mathbb{R})\times J^{k_2}(\mathbb{R},\mathbb{R})\times\cdots\times J^{k_r}(\mathbb{R},\mathbb{R})
\end{equation}
We will study certain actions of finite Lie groups upon such products and show that one can thereby construct an important class of Darboux integrable $n$-hyperbolic manifolds. The main idea can be summarised as follows. Let $\mathbb{J}_1,\mathbb{J}_2$ each be manifolds of the form (\ref{jetProds}). It is natural to consider the distributions on $\mathbb{J}_1$ and $\mathbb{J}_2$ formed by the direct sum of the contact distributions $\mathcal{C}^{k_l}$ of each jet space $J^{k_l}(\mathbb{R},\mathbb{R})$ in a cartesian product $\mathbb{J}$. Let us call these {\it multi-contact distributions} and denote them by symbols $\mathcal{C}(\kappa_1)$ and $\mathcal{C}(\kappa_2)$ on $\mathbb{J}_1$ and $\mathbb{J}_2$, respectively, where
\begin{equation}\label{contact}
\mathcal{C}(\kappa_a):=\mathcal{C}^{k^a_1}\oplus\mathcal{C}^{k^a_2}\oplus\cdots\oplus\mathcal{C}^{k^a_{n_a}},\ a=1,2
\end{equation}
and where the sequences of positive integers
$$
\kappa_a=\langle k^a_1,k^a_2,\ldots,k^a_{n_a}\rangle
$$
completely specify the distributions under consideration. Suppose that free, regular, diagonal actions $\mu_1:G\times \mathbb{J}_1\to \mathbb{J}_1$ and $\mu_2: G\times \mathbb{J}_2\to \mathbb{J}_2$ of a finite Lie group $G$ are symmetries of $(\mathbb{J}_1,\mathcal{C}(\kappa_1))$ and $(\mathbb{J}_2,\mathcal{C}(\kappa_2))$, respectively\footnote{In this paper we ignore those symmetries of $\mathcal{C}(\kappa_a)$ that permute its factors}. Hence, $\mu_1$ and $\mu_2$ restrict to be finite-dimensional contact transformation groups on each factor in $\mathbb{J}_1$ and $\mathbb{J}_2$, respectively. We will confine our study to such actions and show that the quotient of the product $(\mathbb{J}_1\times \mathbb{J}_2,\mathcal{C}(\kappa_1)\oplus\mathcal{C}(\kappa_2))$ by the diagonal action of $G$ determined by $\mu_1$ and $\mu_2$ determines a Darboux integrable $n$-hyperbolic manifold. For more information on symmmetry reduction of exterior differential systems see [\ref{Itskov01}], [\ref{AndersonFels05}], [\ref{AFV09}].  Because all finite dimensional contact transformation groups acting on $J^k(\mathbb{R},\mathbb{R})$ have been classified, [\ref{DoubrovKomrakov}], it is possible to contemplate the problem of classifying all Darboux integrable $n$-hyperbolic manifolds that arise by this proceedure and the PDE systems (\ref{mainEquations}) that they intrinsically determine.
We now describe this construction in detail.

\subsection{The actions on $\mathbb{J}_1$ and $\mathbb{J}_2$}

Express the canonical bases for the multi-contact systems $\mathcal{C}^{(k_i^1)}$ on $\mathbb{J}_1$ by
$$
\Big\{X_i=\P {\alpha_i}+x^i_1\P {x^i}+x^i_2\P {x^i_1}+\cdots+x^i_{k^1_i}\P {x^i_{k^1_i-1}},\ \P {x^i_{k^1_i}}\Big\},\ \ 1\leq i\leq n_1
$$ 
and $\mathcal{C}^{(k_i^2)}$ on $\mathbb{J}_2$ by 
$$
\Big\{Y_j=\P {\beta_j}+y^j_1\P {y^j}+y^j_2\P {x^j_1}+\cdots+y^j_{k^2_j}\P {y^j_{k^2_j-1}},\ \P {y^j_{k^2_j}}\Big\},\ \ 1\leq j\leq n_2.
$$
Using the Fels-Olver method of moving frames [\ref{FO}], we will construct Lie group actions $\mu_1$ on $\mathbb{J}_1$ and $\mu_2$ on $\mathbb{J}_2$ such that 
\begin{enumerate}
\item[A1)] Each action $\mu_1,\mu_2$ is free and regular and $\mu_1$ preserves $\alpha_i$ and $\mu_2$ preserves $\beta_j$ for all $i,j$
\item[A2)] For each $i,\ 1\leq i\leq n_1$ and each $j$,\ $1\leq j\leq n_2$, there is exactly one invariant $u_i$ of $\mu_1$ such that
$$
\frac{\partial u_i}{\partial x^i_{k^1_i}}\neq 0
$$ 
and exactly one invariant $v_j$ such that
$$
\frac{\partial v_j}{\partial y^j_{k^2_j}}\neq 0.
$$ 
\item[A3)] All other invariants of $\mu_a$,\ $a=1,2$,\ and their moving frame components have order lower than 
$\text{min}\,\big(k^a_1,k^a_2,\ldots,k^a_{n_a}\big)$.

\item[A4)] Symmetry group $G$ is chosen so that $\dim\mathbb{J}_1+\dim\mathbb{J}_2-\dim G=3n+1$, where $n=n_1+n_2$.
\end{enumerate}

\begin{defn}
We will call an action of a finite Lie group acting diagonally on a cartesian product $\mathbb{J}$ defined by (\ref{jetProds}) an {\it admissible action on} $\mathbb{J}$ if it is a symmetry of its multi-contact distribution and satisfies conditions A1, A2, A3 and A4.\hfill $\clubsuit$
\end{defn}

\begin{exmp}\label{exC22C2@A1}
Let $G$ be the 2-dimensional non-abelian Lie group. Let $\mathbb{J}_2=J^2(\mathbb{R},\mathbb{R})$ and 
$\mathbb{J}_1=J^2(\mathbb{R},\mathbb{R})\times J^2(\mathbb{R},\mathbb{R})$. The multi-contact systems on $\mathbb{J}_1$ and $\mathbb{J}_2$ are
$$
\begin{aligned}
\mathcal{C}\langle 2,2\rangle&=\{\P {\a_1}+x^1_1\P {x^1}+x^1_2\P {x^1_1},\P {x^1_2}\}\oplus
\{\P {\a_2}+x^2_1\P {x^2}+x^2_2\P {x^2_1},\P {x^2_2}\},\cr
\mathcal{C}\langle 2\rangle &=\{\P {\b_1}+y^1_1\P {y^1}+y^1_2\P {y^1_1},\P {y^1_2}\},
\end{aligned}
$$
respectively. The action of $G$ on $\mathbb{J}_1$ obtained by prolongation to second order of the action on $J^0(\mathbb{R},\mathbb{R})\times J^0(\mathbb{R},\mathbb{R})$ given by
$$
\mu_1(g)(\a_1,x^1,\a_2,x^2)=(\a_1,ax^1+b,\a_2,ax^2+b)
$$
is admissible on $\mathbb{J}_1$, as may be checked. The action of $G$ on $\mathbb{J}_2$ obtained by the prolongation to second order of the action on $J^0(\mathbb{R},\mathbb{R})$ given by
$$
\mu_2(g)(\b_1,y^2)=(\b_1,ay^1+b)
$$
is admissible on $\mathbb{J}_2$. \hfill $\spadesuit$
\end{exmp}

The next step is to construct a {\it right moving frame}, [\ref{FO}], $\rho:\mathbb{J}_1\to G$ for the action $\mu_1$ and a left moving frame
$\lambda:\mathbb{J}_2\to G$, for the action $\mu_2$. 

Let $\a_i,\varpi_1,\varpi_2,\ldots,\varpi_{s_1},u_i$, $1\leq i\leq n_1$ be the differential invariants of $\mu_1$ and let $\b_j, \varsigma_1,\varsigma_2,\allowbreak \ldots,\varsigma_{s_2},v_j$, $1\leq j\leq n_2$ be the differential invariants of $\mu_2$. 
Denote local coordinates on $G$ by $\boldsymbol{w}=w_1,\ldots,w_r$. 
 
The free, regular action $\mu_1$ determines a principal $G$-bundle $\pi^1: \mathbb{J}_1\to \mathbb{J}_1/G:=\mathcal{S}_1$, with fibre diffeomorphic to $G$. Similarly, we have the principal $G$-bundle $\pi^2:\mathbb{J}_2\to \mathbb{J}_2/G:=\mathcal{S}_2$ arising from $\mu^2$. Local trivialisations of $\boldsymbol{\phi^1}$ and $\boldsymbol{\phi^2}$ are local diffeomorphisms
$$
\boldsymbol{\phi^1}:\mathbb{J}_1\to\mathcal{S}_1\times G
$$
and
$$
\boldsymbol{\phi^2}:\mathbb{J}_2\to\mathcal{S}_2\times G
$$
by
\begin{equation}\label{phi1}
\boldsymbol{\phi^1}(\boldsymbol{x})=(\boldsymbol{\a},\boldsymbol{\varpi},\boldsymbol{u},\rho(\boldsymbol{x})),\ \ \ \  \forall\ \ \boldsymbol{x}\in\mathbb{J}_1
\end{equation}
and
\begin{equation}\label{phi2}
\boldsymbol{\phi^2}(\boldsymbol{y})=(\boldsymbol{\b},\boldsymbol{\varsigma},\boldsymbol{v},\lambda(\boldsymbol{y})),
\ \ \ \  \forall\ \ \boldsymbol{y}\in\mathbb{J}_2
\end{equation}
respectively, where,
$$
\begin{aligned}
&\boldsymbol{\alpha}=\a_1,\ldots,\a_{n_1},\ \ \boldsymbol{\beta}=\b_1,\ldots,\b_{n_2},\cr
&\boldsymbol{\varpi}=\varpi_1,\varpi_2,\ldots,\varpi_{s_1},\ \ \boldsymbol{\varsigma}=\varsigma_1,\varsigma_2,\ldots,\varsigma_{s_2},\cr
&\boldsymbol{u}=u_1,\ldots,u_{n_1},\ \ \boldsymbol{v}=v_1,\ldots,v_{n_2}. 
\end{aligned}
$$
and $\boldsymbol{x}$, $\boldsymbol{y}$ denote the canonical multi-contact coordinates on $\mathbb{J}_1$ and  
$\mathbb{J}_2$, respectively. We shall call the maps $\boldsymbol{\phi^1},\boldsymbol{\phi^2}$ the {\it right-} and {\it left augmented moving frames}, respectively.


\begin{defn}
Let $\O$ be a Pfaffian system on smooth manifold $M$ and $\mu : G\times M\to M$ a free, regular Lie group action on $M$ which preserves $\O$,
$$
\mu(g)^*\O\subseteq\O,\ \ \ \  \forall\ \ g\in G.
$$
Then the {\it quotient of} $\O$ {\it by} $G$ is the Pfaffian system $\O/G$ on $M/G$ where,
$$
\O/G=\{\o\in \L^1(M/G)~|~\pi^*\o\in\ \O\},
$$
and where $\pi : M\to M/G$ is the natural projection from $M$ onto the quotient space of $M$ by $G$.\footnote{We assume here that the quotient of a Pfaffian system is also a Pfaffian system; see [\ref{Itskov01}], [\ref{AndersonFels05}] for more details.}\hfill $\clubsuit$
\end{defn}

Naturally, one can likewise define the quotient of a vector field distribution in the obvious way.

\begin{defn}
Let $\mathcal{D}$ be a vector field distribution on manifold $M$ and $\mu : G\times M\to M$ a free and regular Lie group action on $M$ which preserves $\mathcal{D}$,
$$
\mu(g)_*\mathcal{D}\subseteq\mathcal{D},\ \ \ \  \forall\ \ g\in G.
$$
Let $\pi : M\to M/G$ be the natural projection from $M$ onto the quotient space of $M$ by $G$.  Then the {\it quotient of} $\mathcal{D}$ {\it by} $G$ is the distribution on $M/G$ defined by 
$$
\mathcal{D}/G=\ker\,\left(\O/G\right).
$$ 
\end{defn}

Let us denote the images $\boldsymbol{\phi^a}(\mathbb{J}_a)$ by $\S_a$. By construction, for each $a$, $\S_a$ is locally diffeomorphic to $\mathcal{S}_a\times G$, where $\mathcal{S}_1$ has local coordinates 
$(\boldsymbol{\a}, \boldsymbol{\varpi}, \boldsymbol{u})$ 
and $\mathcal{S}_2$ has local coordinates
$(\boldsymbol{\b}, \boldsymbol{\varsigma}, \boldsymbol{v})$. 
It is important to keep distinct the $G$ factors in each of $\S_1$ and $\S_2$. For $\S_1$, we shall denote by $\mathfrak{R}_1$ and $\mathfrak{L}_1$ the infinitesimal right- and left-translations, respectively on its copy of $G$ which we shall sometimes denote $G_{\boldsymbol{a}}$; local coordinates on $G_{\boldsymbol{a}}$ will be denote by $\boldsymbol{a}=(a_1,\ldots,a_r)$. Similarly, infinitesimal right- and left-translations on the $\S_2$ copy of $G$, $G_{\boldsymbol{b}}$, will be denoted by $\mathfrak{R}_2$ and $\mathfrak{L}_2$, respectively and local coordinates on $G_{\boldsymbol{b}}$ will be denote by $\boldsymbol{b}=(b_1,\ldots,b_r)$.

\begin{thm}\label{DMZdistributions1}
For $a=1,2$, let $\mathbb{J}_a$ be products of jet spaces of the form (\ref{jetProds}) with multi-contact distributions $\mathcal{C}(\kappa_a)$.  Suppose admissible actions of a Lie group $G$, $\mu_a:G\times \mathbb{J}_1\times\mathbb{J}_2\to\mathbb{J}_1\times\mathbb{J}_2$ preserve  $\mathcal{C}(\kappa_a)$ for each $a$ and let $\boldsymbol{\phi^a}$ be their right and left augmented moving frames, (\ref{phi1}), (\ref{phi2}). Then, there are functions $\vec{h}_j$,$\vec{k}_j$ on $\mathcal{S}_1$ and functions $\vec{m}_j$,$\vec{n}_j$ on $\mathcal{S}_2$ such that
$$
\begin{aligned}
&\boldsymbol{\phi^1}_*\mathcal{C}(\kappa_1)=\bigoplus_{j=1}^{n_1}\Bigg\{\P {\a_j}+
\vec{h}_j(\boldsymbol{\a},\boldsymbol{\varpi},\boldsymbol{u})\cdot\P {\boldsymbol{\varpi}}+
\vec{k}_j(\boldsymbol{\a},\boldsymbol{\varpi},\boldsymbol{u})\cdot\boldsymbol{L}_1,\ \P {u_j}\Bigg\},\cr
&\boldsymbol{\phi^2}_*\mathcal{C}(\kappa_2)=\bigoplus_{l=1}^{n_2}\Bigg\{\P {\b_l}+
\vec{m}_l(\boldsymbol{\b},\boldsymbol{\varsigma},\boldsymbol{v})\cdot\P {\boldsymbol{\varsigma}}+
\vec{n}_l(\boldsymbol{\b},\boldsymbol{\varsigma},\boldsymbol{v})\cdot\boldsymbol{R}_2,\ \P {v_l}\Bigg\},
\end{aligned}
$$ 
where $\boldsymbol{L}_1$ is a basis for $\mathfrak{L}_1$ and $\boldsymbol{R}_2$ is a basis for $\mathfrak{R}_2$.
\end{thm}
\vskip 3 pt
\noindent{\it Proof.} By suitable prolongations we can arrange that group actions on $\mathbb{J}_a$ are free and regular. By the Olver-Fels method of moving frames [\ref{FO}] we can arrange that the free, regular actions on $\mathbb{J}_1$ and $\mathbb{J}_2$ are admissible. By construction the maps $\boldsymbol{\phi^1}$ and $\boldsymbol{\phi^2}$ are right- and left-equivariant, respectively. That is,
\begin{equation}
\big(\boldsymbol{\phi^1}\circ\mu_1(g)\big)(\boldsymbol{x})=(\boldsymbol{\a},\boldsymbol{\varpi},\boldsymbol{u},
R_{g^{-1}}^1\boldsymbol{w}^1)=\mathbf{R}_{g^{-1}}^1\boldsymbol{\phi^1}(\boldsymbol{x})
\end{equation}        
and
\begin{equation}
\big(\boldsymbol{\phi^2}\circ\mu_2(g)\big)(\boldsymbol{y})=(\boldsymbol{\b},\boldsymbol{\varsigma},\boldsymbol{v},
L_{g}^2\boldsymbol{w}^2)=\mathbf{L}_g^2\boldsymbol{\phi^2}(\boldsymbol{y})
\end{equation}
where $R_g^1,L_g^2$ are right- and left-translations on $G_1$ and $G_2$, respectively and $\mathbf{R}_g^1=\text{id}\times R_g^1$ and 
$\mathbf{L}_g^2=\text{id}\times L_g^2$. We now show that the infinitesimal generators of the action on $\mathbb{J}_1$ pushforward under $\boldsymbol{\phi^1}$ to infinitesimal right-translations $\mathfrak{R}_1$ on $G_1$ while the infinitesimal generators of the action on $\mathbb{J}_2$ pushforward under $\boldsymbol{\phi^2}$ to infinitesimal left-translations $\mathfrak{L}_2$ on $G_2$. 
Also, it will be convenient to let $\boldsymbol{R}_1$, $\boldsymbol{L}_1$ denotes bases for $\mathfrak{R}_1$ and $\mathfrak{L}_1$, respectively.  Similarly $\boldsymbol{R}_2$, $\boldsymbol{L}_2$ denote bases for the Lie algebras $\mathfrak{R}_2$ and $\mathfrak{L}_2$, respectively.

Let $\mathbb{X}_a$ denote bases for the infinitesimal generators of the actions $\mu_a:G\times \mathbb{J}_a\to \mathbb{J}_a$ of $G$ on $\mathbb{J}_a$,\ $a=1,2$. We have
\begin{equation}
\big(\mu_a(g)\big)_*\mathbb{X}_a=\text{Ad}(g)\,\mathbb{X}_a,\ \ \ a=1,2.
\end{equation} 

Hence, in case $a=1$, 
$$
\begin{aligned}
(\mathbf{R}_{g^{-1}}^1)_*\boldsymbol{\phi^1}_*\mathbb{X}_1&=\boldsymbol{\phi^1}_*\mu_1(g)_*\mathbb{X}_1\cr
&=\boldsymbol{\phi^1}_*\left(\text{Ad}(g)\,\mathbb{X}_1\right)\cr
&=\text{Ad}(g)\,\boldsymbol{\phi^1}_*\mathbb{X}_1\cr
&=\left(\mathbf{R}_{g^{-1}}^1\right)_*\left(\mathbf{L}_g^1\right)_*\boldsymbol{\phi^1}_*\mathbb{X}_1,
\end{aligned}
$$
and therefore
$$
\left(\mathbf{L}_g^1\right)_*\boldsymbol{\phi^1}_*\mathbb{X}_1=\boldsymbol{\phi^1}_*\mathbb{X}_1.
$$
Thus $\boldsymbol{\phi^1}_*\mathbb{X}_1$ is a left-invariant $\mathbb{R}^r$-valued vector field on $G_1$. A similar argument shows that $\boldsymbol{\phi^2}_*\mathbb{X}_2$ is a right-invariant $\mathbb{R}^r$-valued vector field on $G_2$.

Now fix any total differential operator $\mathcal{D}_j\in \mathcal{C}(\kappa_1)$ and any infinitesimal generator $X$ of $\mu_1$.
Then since $X$ is an infinitesimal symmetry of $\mathcal{C}(\kappa_1)$, it follows that 
$[X,\mathcal{D}_j]\in \mathcal{C}(\kappa_1)$. By hypothesis, an admissible group action fixes the ``independent variables" $\a_i$ and hence $X\a_i=0$ for all $i$ from which it easily follows that
$$
[X,\mathcal{D}_j]\in\big\{\P {x^1_{k_1}},\ \P {x^2_{k_2}},\ldots,\P {x^{n_1}_{k_{n_1}}}\big\}.
$$
Also, by property A3) for admissible action we deduce that 
$$
\boldsymbol{\phi^1}_*(\P {x_{k^1_j}})=\P {u_j},\ \ \boldsymbol{\phi^2}_*(\P {y_{k^2_l}})=\P {v_l}.
$$
From these facts and because
$\P {\boldsymbol{\a}},\P {\boldsymbol{\varpi}},\P {\boldsymbol{u}},\boldsymbol{L}_1$ form a frame-field on $\boldsymbol{\phi^1}(\mathbb{J}_1)$, we deduce that
$$
\begin{aligned}
\boldsymbol{\phi^1}_*(\mathcal{D}_j)=
\P {\a_j}+\vec{h}_j(\boldsymbol{\a},\boldsymbol{\varpi},\boldsymbol{u},\boldsymbol{a})\cdot\P {\boldsymbol{\varpi}}+
\vec{k}_j(\boldsymbol{\a},\boldsymbol{\varpi},&\boldsymbol{u},\boldsymbol{a})\cdot\boldsymbol{L}_1\cr
&\mod \{\P {u_1},\P {u_2},\ldots,\P {u_n}\}
\end{aligned}
$$
for some functions $\vec{h}_j, \vec{k}_j$ on $\boldsymbol{\phi^1}(\mathbb{J}_1)$. Since $\boldsymbol{\phi^1}_*(X)\in\mathfrak{R}_1$, it follows that
$\vec{h}_j,\vec{k}_j$ are independent of $\boldsymbol{a}$ and so depend only on the invariants of $\mu_1$. They therefore have the form
$$
\boldsymbol{\phi^1}_*(\mathcal{D}_j)=\P {\a_j}+
\vec{h}_j(\boldsymbol{\a},\boldsymbol{\varpi},\boldsymbol{u},)\cdot\P {\boldsymbol{\varpi}}+
\vec{k}_j(\boldsymbol{\a},\boldsymbol{\varpi},\boldsymbol{u},)\cdot\boldsymbol{L}_1\mod \{\P {u_1},\P {u_2},\ldots,\P {u_n}\}.
$$
Similarly, if $\mathcal{E}_l$ is a total differential operator in $\mathcal{C}(\kappa_2)$ then
$$
\boldsymbol{\phi^2}_*(\mathcal{E}_l)=\P {\b_l}+
\vec{m}_l(\boldsymbol{\b},\boldsymbol{\varsigma},\boldsymbol{v})\cdot\P {\boldsymbol{\varsigma}}+
\vec{n}_l(\boldsymbol{\b},\boldsymbol{\varsigma},\boldsymbol{v})\cdot\boldsymbol{R}_2\mod \{\P {v_1},\P {v_2},\ldots,\P {v_n}\}
$$
for some functions $\vec{m}_l,\vec{n}_l$ on $\boldsymbol{\phi^2}(\mathbb{J}_2)$, independent of $\boldsymbol{b}$. 
\hfill\hfill\qed
\vskip 5 pt
Next, let $\mathbf{id}_1 :G_{\boldsymbol{a}}\to G$ and $\mathbf{id}_2: G_{\boldsymbol{b}}\to G$ be identity maps. Let a distribution $H$ be defined on $\mathcal{S}_1\times\mathcal{S}_2\times G$ by
\begin{equation}\label{H}
\begin{aligned}
H=&\bigoplus_{j=1}^{n_1}\Bigg\{\P {\a_j}+
\vec{h}_j(\boldsymbol{\a},\boldsymbol{\varpi},\boldsymbol{u})\cdot\P {\boldsymbol{\varpi}}+
\vec{k}_j(\boldsymbol{\a},\boldsymbol{\varpi},\boldsymbol{u})\cdot\boldsymbol{L},\ \P {u_j}\Bigg\}\oplus\cr
&\bigoplus_{l=1}^{n_2}\Bigg\{\P {\b_l}+
\vec{m}_l(\boldsymbol{\b},\boldsymbol{\varsigma},\boldsymbol{v})\cdot\P {\boldsymbol{\varsigma}}+
\vec{n}_l(\boldsymbol{\b},\boldsymbol{\varsigma},\boldsymbol{v})\cdot\boldsymbol{R},\ \P {v_l}\Bigg\}=H_1\oplus H_2
\end{aligned}
\end{equation}
where $\boldsymbol{L}={\mathbf{id}_1}_*\boldsymbol{L}_1$ and $\boldsymbol{R}={\mathbf{id}_2}_*\boldsymbol{R}_2$ form bases for the infinitesimal left- and right-translations on $G$; for later use we denote their local expression by
$$
L_i=\sum_{j=1}^r\L_i^j(w)\P {w_j},\ \ \ R_i=\sum_{j=1}^r\Pi_i^j(w)\P {w_j},
$$
where $\boldsymbol{w}=(w_1,\ldots,w_r)$ are local coordinates on $G$. 

Next, we define an action on $\Sigma_1\times\Sigma_2=(\mathcal{S}_1\times G_{\boldsymbol{a}})\times (\mathcal{S}_2\times G_{\boldsymbol{b}})$ by
$$
{\mu}_D(g)(\mathcal{I}_1,R_{g^{-1}}\boldsymbol{a},\mathcal{I}_2,L_g\boldsymbol{b}),\ \ \ \forall\ g\in G.
$$
where $\mathcal{I}_1=\boldsymbol{\a},\boldsymbol{\varpi},\boldsymbol{u}$ and $\mathcal{I}_2=\boldsymbol{\b},\boldsymbol{\varsigma},\boldsymbol{v}$. Action $\mu_D$ is a free, regular left-action of $G$ and we denote the quotient of $\S_1\times\S_2$ by the orbits of $\mu_D$ by $\S_1\times_G\S_2$. Observe that 
$$
\pi(\mathcal{I}_1,\boldsymbol{a},\mathcal{I}_2,\boldsymbol{b})=
(\mathcal{I}_1,\mathcal{I}_2,\boldsymbol{a}\cdot\boldsymbol{b})
$$
is $\mu_D$-invariant, mapping $G$-orbits in $\S_1\times\S_2$ to points in $\S_1\times_G\S_2$. Hence 
$$
\pi:\S_1\times\S_2\to\S_1\times_G\S_2
$$ 
is the quotient map.

 We denote the final term in the derived flag of any distribution $\mathcal{E}$ by the symbol $\mathcal{E}^{(\infty)}$. It is the smallest integrable distrubution containing $\mathcal{E}$. By construction, the maximal integral submanifolds of $H_1^{(\infty)}$ are the leaves of the foliation in $\mathcal{S}_1\times\mathcal{S}_2\times G$ defined locally by $\mathcal{I}_2$=constant.  Similarly, those of $H_2^{(\infty)}$ are locally defined by $\mathcal{I}_1$=constant.

Let $S_1$ denote a fixed integral submanifold of $H_2^{(\infty)}$ and $S_2$ a fixed integral submanifold of $H_1^{(\infty)}$ in $\S_1\times_G\S_2\simeq\mathcal{S}_1\times\mathcal{S}_2\times G$. Then as well as $\pi$, define projections $\pi_1,\pi_2$, and inclusions $\iota_1,\iota_2$ such that the following diagram commutes
$$
\begindc{\commdiag}[50]
\obj(1,2){$(\Sigma_1\times\Sigma_2,\Omega_1\oplus\Omega_2)$}
\obj(0,1){$(S_1, \Theta_1)$}
\obj(2,1){$(S_2, \Theta_2$)}
\obj(1,0){$(\Sigma_1\times_G\Sigma_2, \Theta)$}
\mor{$(\Sigma_1\times\Sigma_2,\Omega_1\oplus\Omega_2)$}{$(\Sigma_1\times_G\Sigma_2, \Theta)$}{$\pi$}
\mor{$(\Sigma_1\times\Sigma_2,\Omega_1\oplus\Omega_2)$}{$(S_1, \Theta_1)$}{$\pi_1$}[-1,0]
\mor{$(\Sigma_1\times\Sigma_2,\Omega_1\oplus\Omega_2)$}{$(S_2, \Theta_2$)}{$\pi_2$}
\mor{$(S_1, \Theta_1)$}{$(\Sigma_1\times_G\Sigma_2, \Theta)$}{$\iota_1$}[-1,0]
\mor{$(S_2, \Theta_2$)}{$(\Sigma_1\times_G\Sigma_2, \Theta)$}{$\iota_2$}[1,0]
\enddc
$$
Explicitly
$$
\begin{aligned}
&\pi_1(\mathcal{I}_1,\boldsymbol{a},\mathcal{I}_2,\boldsymbol{\b})=(\mathcal{I}_1, \boldsymbol{a}),\ \ 
\pi_2(\mathcal{I}_1,\boldsymbol{a},\mathcal{I}_2,\boldsymbol{b})=(\mathcal{I}_2, \boldsymbol{b}),\cr
&\iota_1(\mathcal{I}_1,\boldsymbol{w})=(\mathcal{I}_1,\mathcal{I}_2,\boldsymbol{w}),\ \ 
\iota_2(\mathcal{I}_2,\boldsymbol{w})=(\mathcal{I}_1,\mathcal{I}_2,\boldsymbol{w}).
\end{aligned}
$$
We have also defined the following Pfaffian systems
$$
\Theta=H^\perp,\ \Theta_1=\iota_1^*\Theta,\ \ \Theta_2=\iota_2^*\Theta,\ \ \O_1=\pi_1^*\Theta_1,\ \ \O_2=\pi_2^*\Theta_2.
$$

\begin{thm}\label{DMZdistributions2}
Distribution $H$, locally defined by (\ref{H}), can be identified with the quotient of $(\S_1\times\S_2,\O_1\oplus\O_2)$ by $\mu_D$ giving rise to the differential system $(\S_1\times_G\S_2,\Theta)$. The pair
$$
\left(\S_1\times_G \S_2,H\right)
$$ 
is a Darboux integrable $n$-hyperbolic manifold.
\end{thm}
\vskip 3 pt
{\proof 

Pfaffian system $\Theta$ has local expression
\begin{equation}\label{eqQuotient5}
\Big\{\Upsilon_1,\ \Upsilon_2,\ \theta^j\Big\}_{j=1}^r
\end{equation}
where
$$
\theta^j=dw_j-\sum_{l=1}^{r}\sum_{t=1}^{n_1}k^l_t(\mathcal{I}_1)\Lambda^j_l(w)d\alpha_t
-\sum_{l=1}^r\sum_{s=1}^{n_2}n^l_s(\mathcal{I}_2)\Pi^j_l(w)d\beta_s,\ \ 1\leq j\leq r.
$$
The symbol $\Upsilon_1$ denotes 1-forms on $S_1$, while $\Upsilon_2$ denotes 1-forms on $S_2$. 
Denote by $e$ the identity in $G$ and let $(\s,\t)$ be local coordinates around $(e,e)\in G\times G$. In terms of the composition function $\mathbf{c}$ on $G$, we have
\begin{equation}\label{eqQuotient6}
\Lambda^j_i(w)=\frac{\partial \mathbf{c}^j}{\partial {\s_i}}(e,w),\ 
\Pi^j_i(w)=\frac{\partial \mathbf{c}^j}{\partial {\t_i}}(w,e),\ i,j=1,\ldots,r.
\end{equation}
Recalling the right-invariance of the $L_i$ and the left-invariance of the $R_i$, we have, for $\s,\t\in G$
\begin{equation}\label{eqQuotient7}
\begin{aligned}
&(\rho_\t)_*\sum_{j=1}^r\Lambda^j_i(\s)\;{\P {s_j}}_{|_\s}=\sum_{j=1}^r\Lambda^j_i(u\cdot v)\;
{\P {s_j}}_{|_{\s\cdot \t}},\cr 
&(\lambda_\s)_*\sum_{j=1}^r\Pi^j_i(\t)\;{\P {s_j}}_{|_\t}=\sum_{j=1}^r\Pi^j_i(\s\cdot \t)\;{\P {s_j}}_{|_{\s\cdot \t}},
\end{aligned}
\end{equation}
where $\s\cdot \t$ denotes multiplication on $G$ and $\lambda_\s, \rho_\t$ denote left- and right-translations on $G$, respectively. However, we also have
\begin{equation}\label{eqQuotient8}
\begin{aligned}
&(\rho_\t)_*\sum_{j=1}^r\Lambda^j_i(\s)\;{\P {s_j}}_{|_\s}=\sum_{j,k=1}^r\Lambda^k_i(\s)\frac{\partial \mathbf{c}^j}{\partial \s_k}(\s,\t)\;
{\P {s_j}}_{|_{\s\cdot \t}},\cr 
&(\lambda_\s)_*\sum_{j=1}^r\Pi^j_i(\t)\;{\P {s_j}}_{|_\t}=\sum_{j,k=1}^r\Pi^k_i(\t)\frac{\partial \mathbf{c}^j}{\partial \s_k}(\s,\t)\;
{\P {s_j}}_{|_{\s\cdot \t}}.
\end{aligned}
\end{equation}
Recall the projection $\pi :\Sigma_1\times\Sigma_2\to M$ defined by
$$
\pi(\mathcal{I}_1,\boldsymbol{a};\,\mathcal{I}_2,\boldsymbol{b})=
\big(\mathcal{I}_1,\mathcal{I}_2,\mathbf{c}(\boldsymbol{a},\boldsymbol{b})\big).
$$
Upon making use of (\ref{eqQuotient5})-(\ref{eqQuotient8}) we compute that for each $j=1,\ldots,r$
$$
\begin{aligned}
\pi^*\theta^j=&d\mathbf{c}(\boldsymbol{a,b})-\sum_{m,l=1}^r\sum_{t=1}^{n_1}k^l_t(\mathcal{I}_1)\Lambda^m_l(\boldsymbol{a})
\frac{\partial {c}^j}{\partial a_m}\;d\alpha_t\cr
&\hskip 110 pt-\sum_{m,l=1}^r\sum_{s=1}^{n_2}n^l_s(\mathcal{I}_2)\Pi^m_l(\boldsymbol{b})\frac{\partial {c}^j}{\partial v_m}\;d\beta_s\cr
=&\sum_{m=1}^r\frac{\partial {c}^j}{\partial a_m}\Big(
da_m-\sum_{l=1}^r\sum_{t=1}^{n_1}k^l_t(\mathcal{I}_1)\Lambda^m_l(\boldsymbol{a})\;d\alpha_t\Big)\cr
&\hskip 50 pt+\sum_{m=1}^r\frac{\partial {c}^j}{\partial b_m}\Big(db_m-\sum_{l=1}^r\sum_{s=1}^{n_2}n^l_s(\mathcal{I}_2)\Pi^m_l(\boldsymbol{b})\;d\beta_s\Big).
\end{aligned}
$$
On the other hand
$$
\begin{aligned}
&\Omega_1=\pi_1^*\iota_1^*\Theta=\Big\{\Upsilon_1,\ da_j-\sum_{l=1}^r\sum_{t=1}^{n_1}k^l_t(\mathcal{I}_1)\Lambda^j_l(\boldsymbol{a})\;d\alpha_t\Big\}_{j=1}^r\cr
&\Omega_2=\pi_2^*\iota_2^*\Theta=\Big\{\Upsilon_2,\ db_j-\sum_{l=1}^r\sum_{s=1}^{n_2}n^l_s(\mathcal{I}_2)\Pi^j_i(\boldsymbol{b})\;d\beta_s\Big\}_{j=1}^r.
\end{aligned}
$$
We have therefore shown that 
$$
\pi^*\Theta\subseteq\Omega_1\oplus\Omega_2.
$$
That is, $(\S_1\times_G\S_2,\Theta)$ is the quotient of the product structure $(\Sigma_1\times\Sigma_2,\Omega_1\oplus\Omega_2)$ by the action $\mu_D$ as we wanted. 
Thus, we have the identification,
$$
(\S_1\times_G\S_2,\Theta)\approx (\Sigma_1\times\Sigma_2,\Omega_1\oplus\Omega_2 )/G.
$$

We now show that $(\S_1\times_G\S_2,H)$ (equivalently $(\S_1\times_G\S_2,\Theta)$) is a Darboux integrable $n$-hyperbolic manifold. We need to check that $\ch H^{(1)}=\{\P {u_1},\ldots,\P {u_{n_1}},\P {v_1},\ldots,\P {v_{n_2}}\}$ and that the derived type of $H$ is $[[2n,0],[3n,n],[3n+1, 3n+1]]$. Recall that
\begin{equation}
\begin{aligned}
&\ch\mathcal{C}(\kappa_1)^{(1)}=\Big\{\P {x^1_{k_1}}, \P {x^2_{k_2}},\ldots,\P {x^{s}_{k_{s}}}\Big\},\cr
&\ch\mathcal{C}(\kappa_2)^{(1)}=\Big\{\P {y^1_{l_{r_1}}}, \P {y^2_{l_{r_2}}},\ldots,\P {y^{r}_{l_{r}}}\Big\},
\end{aligned}
\end{equation}
where $s=n_1, r=n_2,k=k^1,l=k^2$. Hence,
$$
\ch H_1^{(1)}=\{\P {u_1},\ldots,\P {u_s}\},\ \ \ch H_2^{(1)}=\{\P {v_1},\ldots,\P {v_r}\}.
$$
From the local form of $H$, we see that 
$$
\{\P {u_1},\ldots,\P {u_s}\}\oplus\{\P {v_1},\ldots,\P {v_r}\}\subseteq\ch H^{(1)}.
$$
Let $\chi\in\ch H^{(1)}$ and as $H_1\cap H_2=\{\ 0\ \}$, we have $\chi=\chi_1+\chi_2$ where, for $l=1,2,$ $\chi_l\in H_l$.
But this implies that $\chi_l\in\ch H^{(1)}_l$ and therefore $\chi_l\in(\text{\bf id}_l)_*\boldsymbol{\phi^l}_*\xi_l$ for some
$\xi_l\in\ch\mathcal{C}(\kappa_l)^{(1)}$. Hence $\chi\in\{\P {u_1},\ldots,\P {u_s}\}\oplus\{\P {v_1},\ldots,\P {v_r}\}$. That is,
$$
\ch H^{(1)}\subseteq\{\P {u_1},\ldots,\P {u_s}\}\oplus\{\P {v_1},\ldots,\P {v_r}\}.
$$
Hence, we've shown that 
$\ch H^{(1)}=\{\P {u_1},\ldots,\P {u_s}\}\oplus\{\P {v_1},\ldots,\P {v_r}\}$.

Finally, we must check that $\dim H^{(1)}=3n$ and $\dim H^{(2)}=3n+1$. To establish the first equality, observe that $\dim\mathcal{C}(\kappa_a)^{(1)}=\dim\mathcal{C}(\kappa_a)+n_a$ and therefore $\dim H_a^{(1)}=\dim H_a+n_a$. Now $H^{(1)}_1\cap H^{(1)}_2=\{\ 0\ \}$. For suppose there is an nonzero $Y\in H^{(1)}_1\cap H^{(1)}_2$. This implies, in particular that $Y\in H^{(1)}$ and hence $[Y,H^{(1)}_2]\in H^{(1)}_2$. By the same token, $[Y,H^{(1)}_1]\in H^{(1)}_1$. This implies that 
$Y\in\ch H^{(1)}_1\cap\ch H^{(1)}_2=\{\ 0\ \}$ and hence 
$$
\begin{aligned}
\dim H^{(1)}=\dim H^{(1)}_1+\dim H^{(1)}_2&=\dim H_1+n_1+\dim H_2+n_2\cr
                                          &=2n+n_1+n_2=3n.
\end{aligned}
$$

To see that $\dim H^{(2)}=3n+1$ we note that as $\dim M=3n+1$ and $\dim H^{(1)}=3n$, we have either $\dim H^{(2)}=3n$ or $\dim H^{(2)}=3n+1$. In the former case, $H^{(1)}$ is Frobenius integrable in which case there is a regular function $\eta$ on $M$ which is an invariant of $H^{(1)}$. In particular, $\eta$ is annihilated by each vector field in $H^{(1)}_1$ and each vector field in $H^{(1)}_2$. But the only invariants of $H^{(1)}_1$ are $\boldsymbol{\beta, \varsigma, v}$, while the only invariants of $H^{(1)}_2$ are $\boldsymbol{\alpha, \varpi, u}$.  Hence, there are regular functions $\mathcal{F}_1,\mathcal{F}_2$ such that
$$
\mathcal{F}_1(\boldsymbol{\alpha,\varpi,u})=\eta=\mathcal{F}_2(\boldsymbol{\beta,\varsigma,v}).
$$
This implies that $\mathcal{F}_a$ are constant functions; a contradiction. We have now proved that $(\Sigma_1\times_G\Sigma_2,H)$ is a Darboux integrable $n$-hyperbolic manifold, as we wanted.}
\hfill\hfill\qed
\vskip 5 pt
\begin{exmp} 
This is a continuation of {\it Example \ref{exC22C2@A1}}. To simplify notation we rewrite the multi-contact systems as
$$
\begin{aligned}
\mathcal{C}\langle 2,2\rangle &=\{\P x+x_1\P {x_0}+x_2\P {x_1},\P {x_2}\}\oplus
\{\P y+y_1\P {y_0}+y_2\P {y_1},\P {y_2}\}\cr
\mathcal{C}\langle 2\rangle &=\{\P z+z_1\P {z_0}+z_2\P {z_1},\P {z_2}\}
\end{aligned}
$$
Then $\mu_1(g):\mathbb{J}_1\to \mathbb{J}_1$ is given by
$$
\mu_1(a,b)(\boldsymbol{x},\boldsymbol{y})=(x,ax_0+b,ax_1,ax_2,y,ay_0+b,ay_1,ay_2)
$$
and $\mu_2(g):\mathbb{J}_2\to \mathbb{J}_2$ is given by
$$
\mu_2(a,b)(\boldsymbol{z})=(z,az_0+b,az_1,az_2).
$$
A complete list of $\mu_1$ invariants is
$$
x,y,\frac{x_0-y_0}{x_1},\frac{y_1}{x_1},\frac{x_2}{x_1},\frac{y_2}{y_1}
$$
whose right-moving frame is
$$
\rho : \mathbb{J}_1\to G
$$
given by
$$
\rho(\boldsymbol{x},\boldsymbol{y})=\left(\begin{matrix}\frac{1}{x_1}\ &\ -\frac{x_0}{x_1}\cr
0 & 1\end{matrix}\right)=\left(\begin{matrix}\  a\  &\  \ b\ \cr
                                   0\ &\ 1\end{matrix}\right).
                                   $$
A complete list of $\mu_2$ invariants is
$$
z, \frac{z_2}{z_1}
$$                                   
whose left moving frame is                                   
$$
\lambda(\boldsymbol{z})=\left(\begin{matrix}\frac{1}{z_1}\ &\ -\frac{z_0}{z_1}\cr
0 & 1\end{matrix}\right)^{-1}=\left(\begin{matrix} z_1\ &\ -z_0\cr
                                   0\ &\ 1\end{matrix}\right).
                                   $$
Consequently, the diffeomorphism $\boldsymbol{\phi^1}:\mathbb{J}_1\to\mathbb{R}^6\times G_{\boldsymbol{a}}$ for this example is
$$
\boldsymbol{\phi^1}(\boldsymbol{x},\boldsymbol{y})=\left(x,y,\frac{x_0-y_0}{x_1},\frac{y_1}{x_1},\frac{x_2}{x_1},\frac{y_2}{y_1},
\frac{1}{x_1}, -\frac{x_0}{x_1}\right)=(x,y,\varpi_1,\varpi_2,u_1,u_2,a_1,a_2)
$$
while diffeomorphism $\boldsymbol{\phi^2} : \mathbb{J}_2\to \mathbb{R}^2\times G$ is
$$
\boldsymbol{\phi^2}(\boldsymbol{z})=\left(z,\frac{z_2}{z_1},z_1,-z_0\right)=(z,v_1,b_1,b_2).
$$
The pushforward of $\mathcal{C}(2,2)$ by $\boldsymbol{\phi^1}$ yields
$$
\begin{aligned}
\boldsymbol{\phi^1}_*\mathcal{C}\langle 2,2\rangle=&\big\{\P x+(1-\varpi_1u_1)\P {\varpi_1}-\varpi_2u_1\P {\varpi_2}-u_1^2\P {u_1}-u_1a_1\P {a_1}-(1+u_1a_2)\P {a_2}, \P {u_1}\big\}\cr
&\hskip 140 pt\oplus\{\P y-\varpi_2\P {\varpi_1}+\varpi_2u_2\P {\varpi_2}-u_2^2\P {u_2},\P {u_2}\}\cr
=&\{\tilde{X}_1,\P {u_1}\}\oplus\{\tilde{X}_2,\P {u_2}\},
\end{aligned}
$$
and the pushforward of $\mathcal{C}\langle 2\rangle$ by $\boldsymbol{\phi^2}$ yields
$$
\begin{aligned}
\boldsymbol{\phi^2}_*\mathcal{C}\langle 2\rangle=&\big\{\P z-v_1^2\P {v_1}+vb_1\P {b_1}-b_1\P {b_2},\P {v_1}\big\}\cr
                      =&\{\tilde{Y}_1,\P {v_1}\}.
\end{aligned}
$$
Next, in accordance with Theorem \ref{DMZdistributions2},  we form the differential system
$(\S_1\times_G \S_2,H)$ with
$$
H=\{X_1,\P {u_1}\}\oplus\{X_2,\P {u_2}\}\oplus\{Y_1,\P {v_1}\}
$$
where
$$
\begin{aligned}
&X_1=\P x+(1-\varpi_1u_1)\P {\varpi_1}-\varpi_2u_1\P {\varpi_2}-u_1^2\P {u_1}-u_1w_1\P {w_1}-(1+u_1w_2)\P {w_2},\cr
&X_2=\P y-\varpi_2\P {\varpi_1}+\varpi_2u_2\P {\varpi_2},\cr
&Y_1=\P z-v_1^2\P {v_1}+vw_1\P {w_1}-w_1\P {w_2}.
\end{aligned}
$$
It can be checked that $(\S_1\times_G \S_2,H)$ is a Darboux integrable 3-hyperbolic manifold. Indeed, its derived type is
$$
[[6,0],[9,3],[10,10]]
$$
and
$$
\ch H^{(1)}=\{\P {u_1},\P {u_2},\P {v_1}\},
$$
proving that $H$ is a Darboux integrable 3-hyperbolic distribution, since it also has the correct local normal form. In accordance with Theorem \ref{DMZdistributions2}, $H$ is the quotient of 
$\mathcal{C}\langle 2,2\rangle\oplus\mathcal{C}\langle 2\rangle$ by the diagonal action and the natural map
$$
\pi:\S_1\times \S_2\to \S_1\times_G \S_2
$$
is given by
$$
\pi(\boldsymbol{x},\boldsymbol{y},\boldsymbol{z})=(x,y,\varpi_1,\varpi_2,u_1,u_2,z,v_1,g_1,g_2)
$$
where $g_1,g_2$ satisfy
$$
\left(\begin{matrix}\, g_1\ &\ g_2\,\cr
                      0\ &\ 1\end{matrix}\right)=\left(\begin{matrix} a_1\ &\ a_2\cr
                                                                       0  &  1\end{matrix}\right)
                                                                       \left(\begin{matrix} b_1\ &\ b_2\cr
                                                                       0  &  1\end{matrix}\right).
$$

We can now implement Theorem \ref{KTdistributions}. We find that
$$
\mathcal{A}=\{\P {u_1},\P {u_2},\P {v_1},\P {\varpi_2},\P {w_1},\varpi_1\P {\varpi_1}+w_2\P {w_2}\}
$$
whose invariants are $x,y,z,p=\varpi_1/w_2$. Computing coordinates as in Theorem \ref{KTdistributions} and relabelling $Y_1$ as $X_3$
$$
p_i=X_ip,p_{ij}=X_ip_j
$$
we find
\begin{equation}\label{Ex3.1defingEqs1}
p=\varpi_1/w_2,p_1=\varpi_1/w_2^2+1/w_2,p_2=-\varpi_2/w_2,p_3=\varpi_1w_1/w_2^2
\end{equation}
and
\begin{equation}\label{Ex3.1defingEqs2}
p_{12}=f_{12}=-\varpi_2/w_2^2,p_{13}=f_{13}=w_1(w_2+2\varpi_1)/w_2^3,p_{23}=f_{23}=-\varpi_2w_1/w_2^2
\end{equation}
Solving equation (\ref{Ex3.1defingEqs1}) for $\varpi_1,\varpi_2,w_1,w_2$ in terms of $p,p_1,p_2,p_3$ and substituting into equation (\ref{Ex3.1defingEqs2}) reveals the {\sf GDMZ} system in standard jet coordinates
\begin{equation}\label{Ex3.1PDE}
\boxed{
u_{xy}=\frac{2u+1}{u(u+1)}u_xu_z,\ u_{xz}=\frac{1}{u+1}u_xu_y,\ u_{yz}=\frac{1}{u}u_yu_z
}
\end{equation}
As predicted by Theorem \ref{KTdistributions}, system (\ref{Ex3.1PDE}) is semilinear and involutive.\hfill $\spadesuit$
\end{exmp}

By its construction, the {\sf GDMZ} system (\ref{Ex3.1PDE}) is {\it Darboux integrable}. Using the methods of [\ref{AFV09}], its general solution can be expressed in finite terms of three arbitrary functions, each of one variable. This fact will be used in the next section when solutions of a given {\sf DMZ} system are required.

Though the characterisation of {\sf GDMZ} systems provided by Theorem \ref{KTdistributions} is quite general our specific interest in this paper is in {\it linear} {\sf GDMZ} systems since these are the ones that are known to have interesting applications. To construct these, we use the earlier method of quotienting products of jet spaces by a group action but now choosing the $G$-action to be that of an {\it abelian} group. To ilustrate the significance of this, we now briefly indicate how to use it to construct solutions of {\sf 3WRI} and {\sf m3WRI} systems. A complete treatment is presented in [\ref{SergeevVassiliou09}]

\section{Solutions of {\sf 3WRI} systems from 3-hyperbolic manifolds. Role of gauge transformations}
\label{sec:3}

The purpose of this section is to give an indication of how the theory developed so far can be applied to construct solutions for the $n$-wave resonant interaction system and the modified $n$-wave resonant interaction system that we discussed briefly in the Introduction. For simplicity of exposition, the case $n=3$ will be discussed but the reader will quickly surmise its generalisation to arbitrary $n$.
 
We consider the following action of $\mathbb{R}^4$ on $J^4(\mathbb{R},\mathbb{R})$
$$
\begin{aligned}
\mu^1_g(z,q,q_1,q_2,q_3,q_4)\mapsto \Big(z,q-\frac{z^2}{2}t_1+\frac{z}{2}t_2-&\frac{t_3}{6}+\frac{z^3}{6}t_4,\cr 
                                  q_1-zt_1+\frac{1}{2}t_2+\frac{z^2}{2}t_4,q_2-&t_1+zt_4,q_3+t_4,q_4\Big),
\end{aligned}
$$
where we have deviated a little from the notation of section 3 by denoting the standard coordinates on $J^4(\mathbb{R},\mathbb{R})$ by $z,q,q_1,q_2,q_3,q_4$ and coordinates on $G=\mathbb{R}^4$ are $t_1,\ldots,t_4$.

Denote the standard coordinates on $J^2(\mathbb{R},\mathbb{R})\times J^2(\mathbb{R},\mathbb{R})$ by $x,m,m_1,m_2,y,n,n_1,n_2$ and define the action of $\mathbb{R}^4$  by
$$
\begin{aligned}
\mu^2_g(x,m,m_1,m_2,y,n,n_1,n_2)=(x,m+t_2-xt_1,&m_1-t_1,m_2,\cr
                                    &y,n+t_4-yt_3,n_1-t_3,n_2)
\end{aligned}
$$
The invariants of $\mu^1$ are $z,q_2$ while the $G$-equivariant moving frame is
$$
\rho(z,q,q_1,q_2,q_3,q_4)=(q_2-zq_3,2zq_2-z^2q_3-2q_1,6q-6zq_1+3z^2q_2-z^3q_3,-q_3).
$$
The invariants of $\mu_2$ are $x,m_2,y,n_2$ and the $G$-equivariant moving frame is
$$
\lambda(x,m,m_1,m_2,y,n,n_1,n_2)=(m_1,xm_1-m,n_1,yn_1-n).
$$
Forming the diffeomorphisms $\boldsymbol{\phi^1,\phi^2}$ are indicated we find that
$$
\boldsymbol{\phi^1}_*\mathcal{C}\langle 4\rangle=\Big\{\P z-u(z\P {a_1}+z^2\P {a_2}+z^3\P {a_3}+\P {a_4}),\ \P u\Big\}
$$
and
$$
\boldsymbol{\phi^2}_*\mathcal{C}\langle 2,2\rangle=\Big\{\P x+v_1(\P {b_1}+x\P {b_2}),\ \P y+v_2(\P {b_3}+y\P {b_4}),\ \P {v_1},\ \P {v_2}\}
$$

As proved in section 3, the differential system $(M,H)$ where 
$$
H=\{X,\P {v_1}\}\oplus\{Y,\P {v_2}\}\oplus\{Z,\P u\}
$$ 
on
$$
M=\boldsymbol{\phi^1}\Big(J^2(\mathbb{R},\mathbb{R})\times J^2(\mathbb{R},\mathbb{R})\Big)\times_G\boldsymbol{\phi^2}\Big(J^2(\mathbb{R},\mathbb{R})\Big)
$$
is a Darboux integrable 3-hyperbolic manifold, where
$$
\begin{aligned}
&X=\P x+v_1(\P {w_1}+x\P {w_2}),\ \ Y=\P y+v_2(\P {w_3}+y\P {w_4})\cr
&Z=\P z-u(z\P {w_1}+z^2\P {w_2}+z^3\P {w_3}+\P {w_4}).
\end{aligned}
$$
Implementing the procedure in Theorem \ref{KTdistributions} we compute the integrable distribution
$$
\mathcal{A}=\{\P {v_1}, \P {v_2}, \P u, \P {w_1}+x\P {w_2}, \P {w_3}+y\P {w_4}, z\P {w_1}+z^2\P {w_2}+z^3\P {w_3}+\P {w_4}\}.
$$
Its invariants are spanned by 
$$
x,\ y,\ z,\ (yz^3-1)(xw_1-w_2)+z(z-x)(yw_3-w_4)
$$ 
According to Theorem \ref{KTdistributions} we can take 
$$
u=(yz^3-1)(w_2-xw_1)+z(z-x)(w_4-yw_3)
$$ 
as the dependent variable after which differentiation by $x-,\, y-,\, z-\,$ total differential operators
$X, Y, Z$,  respectively generates all the  higher order jet variables 
$$
\begin{aligned}
&u_x=-yw_1z^3-zw_4+zyw_3+w_1,\ \ u_y=z^3(w_2-xw_1)-z(z-x)(1+z)w_3,\cr
&u_z=3yz^2(w_2-xw_1)+2zw_4-2zyw_3-xw_4+xyw_3,\cr
\cr
&u_{xy}=z(w_3-z^3w_1),\ \ u_{xz}=-3yz^2w_1-w_4+yw_3,\cr
&u_{yz}=3z^2w_2-3z^2xw_1-2w_3z+xw_3,
\end{aligned}
$$
and leads to the linear {\sf DMZ} system
\begin{equation}\label{m3WRIsol}
\begin{aligned}
&u_{xy}+\frac{z^3}{1-z^3y}u_x-\frac{1}{x-z}u_y-\frac{z^3}{(x-z)(1-z^3y)}u=0,\cr
&u_{xz}+\frac{3yz^2}{1-z^3y}u_x-\frac{1+2yz^3}{2xyz^3+x-2z-yz^4}u_z-\cr
&\hskip 150 pt\frac{3yz^2(1+2yz^3)}{(1-z^3y)(2xyz^3+x-2z-yz^4)}u=0\cr
&u_{yz}-\frac{x-2z}{z(x-z)}u_y-\frac{z^3(2x-z)}{2xyz^3+x-2z-yz^4}u_z+\cr
&\hskip 150 pt \frac{(x-2z)(2x-z)z^2}{(x-z)(2xyz^3+x-2z-yz^4)}u=0.
\end{aligned}
\end{equation}
It can be check, in accordance with previous results that this overdetermined system of three equations in one unknown is a {\sf DMZ} system. In particular, its coefficients satisfy the integrability conditions (\ref{DMZfullCompatibility}) and is therefore involutive.

Note that system (\ref{m3WRIsol}) satisfies the constraint 
\begin{equation}\label{SergeevConstraint}
C_{ij}=\G_{ij}\G_{ji},\ \  \text{for}\ \ \ (i,j)=(1,2),\ (1,3),\, (2,3),
\end{equation}
and therefore the coefficients $\G_{ij}$ in (\ref{m3WRIsol}) constitute a solution of the {\it modified} {\sf 3WRI} system.  From this example, it is evident that a great many {\sf DMZ} systems can be constructed by applying the results of sections 2 and 3.

\subsection{Gauge transformations} Even if the products of jet spaces and the diagonal action is fixed, there is still freedom in the choice of invariant $u$ that ultimately becomes the dependent variable. This is because we have the freedom to replace $u$ by an arbitrary function of $x,y,z$ times $u$: $u\mapsto \lambda(x,y,z)u$, while preserving the linearity of the resulting {\sf GDMZ} system. This ``gauge transformation" raises the important question as to how to construct those special {\sf DMZ} systems that lead to solutions of the {\sf 3WRI} system or the {\sf m3WRI} system.  

To this end, recall that the general Darboux-Manakov-Zakharov linear problem is the {\sf DMZ} system
$$
\left(\P {x_i}\P {x_j}-\G_{ij}\P {x_j}-\G_{j\,i}\P {x_i}+
C_{ij}\right)\phi=0,\ \ 1\leq i<j\leq 3,
$$
where $\G_{ij},\G_{j\,i},C_{ij}$ are functions of $x_1,x_2,x_3$. The case $C_{ij}=0$ is the linear problem for the {\sf 3WRI} system while the case $C_{ij}=\G_{ij}\G_{j\,i}$ is the linear problem for the {\sf m3WRI} system. How can we construct {\sf DMZ} systems that fall into these special classes? To answer this question it is natural to invoke the study of gauge transformations and gauge equivalence of differential operators. Thus, we define the matrix-valued differential operator
\begin{equation}\label{involutiveOperator}
\mathcal{D}=\left(\begin{matrix} 
\Q {x_1,x_2}-\G_{21}\P {x_1}-\G_{12}\P {x_2}+C_{12}\cr
\Q {x_1,x_3}-\G_{31}\P {x_1}-\G_{13}\P {x_3}+C_{13}\cr
\Q {x_2,x_3}-\G_{32}\P {x_2}-\G_{23}\P {x_3}+C_{23}
\end{matrix}\right)=
\left(\begin{matrix} L_{12}\cr
                     L_{13}\cr
                     L_{23}\end{matrix}\right)
\end{equation}
so that $\mathcal{D}\,u=0$ is a {\sf DMZ} system. We shall call such an operator $\mathcal{D}$ {\it involutive}. Each operator element $L_{ij}$ 
can be expressed as a composition of first order operators modulo a correction term
$$
L_{ij}=\big(\P {x_i}-\G_{j\,i}\big)\circ\big(\P {x_j}-\G_{ij}\big)+h_{ij}(\mathcal{D}),\ \ \ (i,j)\in\text{perm}_2\{1,2,3\} 
$$
where
$$
h_{ij}(\mathcal{D})=\frac{\pl \G_{ij}}{\pl x_j}-\G_{ij}\G_{ji}+C_{ij}.
$$
We define the infinite Lie pseudogroup of gauge transformations acting on operators $\mathcal{D}$ by
$$
\bar{\mathcal{D}}=\mathbf{T}_\lambda\mathcal{D}=e^{-\lambda}\mathcal{D}\circ e^\lambda=\left(\begin{matrix} 
\Q {x_1,x_2}-\bar{\G}_{21}\P {x_1}-\bar{\G}_{12}\P {x_2}+\bar{C}_{12}\cr
\Q {x_1,x_3}-\bar{\G}_{31}\P {x_1}-\bar{\G}_{13}\P {x_3}+\bar{C}_{13}\cr
\Q {x_2,x_3}-\bar{\G}_{32}\P {x_2}-\bar{\G}_{23}\P {x_3}+\bar{C}_{23}
\end{matrix}\right)
$$
where $\lambda(x)$ is an arbitrary function of $x_i$ in the domain of the operator $\mathcal{D}$. The coefficients of the image operator are
\begin{equation}\label{gaugedCoefficients}
\bar{\G}_{ij}=\G_{ij}-\P {x_j}\lambda,\ \ \ \ \ \ \bar{C}_{(ij)}=L_{ij}\lambda.
\end{equation}

It is straightforward to prove that the functions $h_{ij}: \boldsymbol{\mathcal{Z}}\to C^{\infty}(\mathbb{R}^3)$ viewed as maps from the set of differential operators $\boldsymbol{\mathcal{Z}}$ of the form $\mathcal{D}$, such that $\mathcal{D}u=0$ is a {\sf DMZ} system,  to the smooth real-valued functions $C^{\infty}(\mathbb{R}^3)$ on $\mathbb{R}^3$ are invariants of the gauge action
$$
h_{ij}(\mathbf{T}_\lambda\mathcal{D})=h_{ij}(\mathcal{D}), \ \ \forall\ \lambda\in C^{\infty}(\mathbb{R}^3),\ \ \text{and}\ \ \mathcal{D}\in\boldsymbol{\mathcal{Z}}.
$$ 
In fact functions $h_{ij}$ are {\it complete invariants} in the sense that operators $\mathcal{D}$ and $\mathcal{D}'$ are gauge equivalent, $\mathcal{D}'=\mathbf{T}_\lambda\mathcal{D}$ for some $\lambda$ if and only if the all their gauge invariants agree
$$
h_{ij}(\mathcal{D}')=h_{ij}(\mathcal{D}),\ \ \forall\ \ (i,j)\in\text{perm}_2\{1,2,3\}.
$$

\begin{prop}
If an operator $\mathcal{D}$ is involutive then so is each operator in its gauge orbit.
\end{prop}
\proof Gauge equivalence of operators is derived from the prolongation to second order of the local diffeomorphism 
$$
\phi: x_1\mapsto x_1,\ x_2\mapsto x_2,\ x_3\mapsto x_3,\ u\mapsto \lambda(x_1,x_2,x_3)u,
$$ 
a contact equivalence of the differential equations that they define. Clearly $\phi$ preserves the independence form $dx_1\wedge dx_2\wedge dx_3$. The image of an involutive differential equation under a contact transformation which preserves the independence form is also involutive with the same independence form. \hfill\qed 
\vskip 5 pt
A general involutive operator (\ref{involutiveOperator}) generates a solution of the {\sf 3WRI} equations once a solution of $\mathcal{D}u=0$ is known. Because the involutive operators arising from our construction are Darboux integrable, it follows that solutions of $\mathcal{D}u=0$ can be explicitly constructed.
\begin{prop}\label{General2_3wave}
Let $\mathcal{D}$ be any involutive differential operator (\ref{involutiveOperator}). Let $u$ be any nonzero solution of $\mathcal{D}u=0$. Then the coefficients of $\mathcal{D}_0=\mathbf{T}_u\mathcal{D}$ provide a solution of the {\sf 3WRI} system.
\end{prop}
\proof Formulas (\ref{gaugedCoefficients}) and the previous Proposition prove
that $\mathcal{D}_0$ is an involutive operator with $\bar{C}_{ij}=0$. Theorem \ref{geometricBackgroundThm} proves that the remaining nonzero coefficients of
$\mathcal{D}_0$ provide a solution of the {\sf 3WRI} system via the Lame potentials $h_1,h_2,h_3$.\hfill\qed  
\vskip 5 pt
Next, we point out that an involutive operator (\ref{involutiveOperator}) satisfying $C_{ij}=0$ for all $i,\ j$, determines a solution of the {\it modified} {\sf 3WRI} system, that is a solution of (\ref{m3WRI}) once a certain quadrature can be performed.
\vskip 5 pt 
\begin{prop}
Suppose $\mathcal{D}_0\in\mathcal{Z}$ is an involutive operator (\ref{involutiveOperator}) with $C_{ij}=0$ for all $(i,j)$. Then 
\begin{enumerate}
\item The linear inhomogeneous PDE system
\begin{equation}\label{SergeevGauge}
\frac{\pl^2 \lambda}{\pl x_i\pl x_j}=\G_{ij}\G_{j\,i},\ \ \  (i,j)\in\{(1,2),\ (1,3),\ (2,3)\},
\end{equation}
for function $\lambda(x_1,x_2,x_3)$ is involutive, and 
\item If $\lambda$ satisfies (\ref{SergeevGauge}) then the coefficients of $T_\lambda\mathcal{D}_0$ solve the modified {\sf 3WRI} system.
\end{enumerate}
\end{prop}
\proof In order for $\mathcal{D}_0$ to be transformed by some gauge to the linear problem $\mathcal{D}_1$ for the {\sf m3WRI} system we require
$$
\mathcal{D}_0\lambda=\left(\begin{matrix} (\G_{21}-\lambda_{x_2})(\G_{12}-\lambda_{x_1})\cr
                                          (\G_{31}-\lambda_{x_3})(\G_{13}-\lambda_{x_1})\cr
                                          (\G_{32}-\lambda_{x_3})(\G_{23}-\lambda_{x_2})\end{matrix}\right)
                     =\left(\begin{matrix} h_{21}(\mathbf{T}_\lambda\mathcal{D}_0)-\P {x_1}(\G_{21}-\lambda_{x_2})\cr
                                           h_{31}(\mathbf{T}_\lambda\mathcal{D}_0)-\P {x_1}(\G_{31}-\lambda_{x_3})\cr
                                           h_{32}(\mathbf{T}_\lambda\mathcal{D}_0)-\P {x_2}(\G_{21}-\lambda_{x_3})\end{matrix}\right).
$$
Since the $h_{ij}$ are gauge-invariant this simplifies to equation (\ref{SergeevGauge})
whose integrability conditions are
\begin{equation}\label{mDMZtransIC}
\frac{\pl}{\pl x_k}\left(\G_{ij}\G_{ji}\right)=\frac{\pl}{\pl x_j}\left(\G_{ik}\G_{ki}\right),\ \ \ (i,j,k)\in\ \text{perm}(1,2,3).
\end{equation}
We recall that the $\G_{ij}$ satisfy the integrability conditions (\ref{KT_0Compatibility}). 
An calculation shows that the integrability conditions (\ref{mDMZtransIC}) are satisfied modulo (\ref{KT_0Compatibility}).\hfill\qed                                      
\vskip 5 pt
Thus each solution of the {\sf 3WRI} system gives rise to solutions of the {\sf m3WRI} system.


\begin{exmp}
Proposition \ref{General2_3wave} asserts that every involutive operator determines another, $\mathcal{D}_0$, with $C_{ij}=0$.  In turn $\mathcal{D}_0$ determines a solution of the 3-wave resonant interaction system as described in Theorem \ref{geometricBackgroundThm}. We shall now illustrate this construction beginning with the involutive system (\ref{m3WRIsol}) which was constructed by the theory of sections 2 and 3.
Performing the gauge transformation $\mathbf{T}_\lambda\mathcal{D}$, with 
$\lambda=-\ln {z(x-z)}$ produces the operator
\begin{equation}\label{newDMZ}
\bar{\mathcal{D}}=\mathbf{T}_\lambda\mathcal{D}=\left(\begin{matrix} 
\Q {x_1,x_2}-\bar{\G}_{21}\P {x_1}-\bar{\G}_{12}\P {x_2}\cr
\Q {x_1,x_3}-\bar{\G}_{31}\P {x_1}-\bar{\G}_{13}\P {x_3}\cr
\Q {x_2,x_3}-\bar{\G}_{32}\P {x_2}-\bar{\G}_{23}\P {x_3}
\end{matrix}\right)
\end{equation}
where
$$
\begin{aligned}
&\bar{\G}_{21}=\frac{z^3}{yz^3-1},\cr 
&\bar{\G}_{31}=\frac{x+2xyz^3-yz^4-2z}{z(x-z)(yz^3-1)},\ 
\ \bar{\G}_{13}=\frac{z(1-yz^3)}{(x+2xyz^3-yz^4-2z)(x-z)},\cr 
&\phantom{\bar{\G}_{31}=\frac{x+2xyz^3-yz^4-2z}{z(x-z)(yz^3-1)}}\ \ \ \ \bar{\G}_{23}=\frac{z^3(2x-z)}{x+2xyz^3-yz^4-2z},
\end{aligned}
$$
with $\bar{\G}_{12}=\bar{\G}_{32}=0$. From these coefficients we obtain Lame potentials 
$$
h_1=\frac{yz^3-1}{z(z-x)},\ \ h_2=1,\ \ h_3=\frac{-z^4y-2z+x+2xz^3y}{x-z},
$$
and the corresponding solution of the $2+1$\,-\,dimensional 3-wave resonant interaction equations, in accordance with Theorem \ref{geometricBackgroundThm} is easily computed to be
\begin{equation}\label{3WRIsol0}
\boxed{
\left(\begin{matrix}
0&A_{12}&A_{13}\cr
A_{21}&0&A_{23}\cr
A_{31}&A_{32}&0  
\end{matrix}\right)
=\frac{1}{x-z}\left(\begin{matrix} 0&0 & z^2\cr 
-z^2y^{-1} & 0 & z^3(2x-z)y^{-1}\cr 
-z^{-2} &0 & 0\end{matrix}\right)}
\end{equation}

Because of Proposition \ref{General2_3wave}, we can generate new solutions of the 2+1-\newline dimensional {\sf 3WRI} system from any given solution arising from an appropriate operator $\mathcal{D}_0$, such as we have constructed in this example, by computing the image of $\mathcal{D}_0$ under a gauge transformation, $\mathbf{T}_u\mathcal{D}_0$, where $u$ is any nonconstant solution of $\mathcal{D}_0u=0$. Normally there is little hope of finding any such solutions but the situation is quite otherwise here because the operator $\mathcal{D}_0$ arises from a Darboux integrable $n$-hyperbolic distribution (with $n=3$ in this case) 
$$
\begin{aligned}
H=\Big\{\P x+v_1(\P {w_1}+x\P {w_2}),\ \ \P {v_1}\Big\}\oplus&\Big\{\P y+v_2(\P {w_3}+y\P {w_4}),\P a,\ \ \P {v_2}\Big\}\oplus\cr
\Big\{\P z-u&(z\P {w_1}+z^2\P {w_2}+z^3\P {w_3}+\P {w_4}),\ \ \P u\Big\}\cr
&\hskip 45 pt =H_1\oplus H_2\oplus H_3.
\end{aligned}
$$
This distribution can be integrated using methods of [\ref{AFV09}]. By construction its {\it Vessiot group} [\ref{AFV09}] is $\mathbb{R}^4$. Furthermore $H_1,H_2$ and $H_3$ are each easily integrated and the fact that the Vessiot group is $\mathbb{R}^4$ means that the superposition formula is the usual linear superposition. From this we easily obtain the general integral submanifold of $H$ to be
$$
\begin{aligned}
\sigma(x,y,z)&=\cr
&\Big(-A_{xx},\ -B_{yy},\ -\frac{1}{6}z^3C_{zzzz}-2z^2C_{zzz}-6C_{zz}z-4C_{z},\cr
&\frac{1}{6}z^4C_{zzz}+\frac{4}{3}z^3C_{zz}+2z^2C_z-A_x,\cr 
&\frac{7}{6}z^4C_{zz}+\frac{1}{6}z^5C_{zzz}+\frac{4}{3}z^3C_z+A-xA_x,\cr 
&z^4C_z+\frac{1}{6}z^6C_{zzz}+z^5C_{zz}-B_y,\cr 
&\frac{1}{6}z^3C_{zzz}+\frac{3}{2}z^2C_{zz}+3C_zz+C+B-yB_y\Big)\cr
&=(a,b,c,w_1,w_2,w_3,w_4).
\end{aligned}
$$ 
where $A(x),\ B(y),\ C(z)$ are arbitrary real-valued functions. We recall that 
$$
u=\frac{(yz^3-1)(xw_1-w_2)+z(z-x)(yw_3-w_4)}{z(x-z)}
$$
generates the involutive operator (\ref{newDMZ}) and hence 
\begin{equation}\label{DMZsol}
\begin{aligned}
\lambda=&\sigma^*u=\cr
&\Big(6xyz^5C_z-2z^4C_{zz}-10z^3C_z-6z^2C-6z^2B+2xyz^6C_{zz}+xz^3C_{zz}+6xz^2C_z-\cr
           &yz^7C_{zz}-2yz^6C_z-6yz^3A+6xzC+6xzB+6A\Big)\Big(z(x-z)\Big)^{-1}
\end{aligned}
\end{equation}
is the general solution of $\bar{\mathcal{D}}\lambda=0$. It follows that the coefficients of $e^{-\lambda}\bar{\mathcal{D}}\circ e^\lambda$ determine a new solution of the 2+1-dimensional {\sf 3WRI} equations by Proposition \ref{General2_3wave}.

For instance, choosing $A(x)=B(y)=C(z)=1$ yields the particular solution
$$
\lambda=\frac{6-12z^2-6yz^3+12xz}{z(x-z)}
$$
of $\bar{\mathcal{D}}\lambda=0$. We compute
\begin{equation}
\widehat{\mathcal{D}}=\mathbf{T}_\lambda\bar{\mathcal{D}}=\left(\begin{matrix} 
\Q {x_1,x_2}-\hat{\G}_{21}\P {x_1}-\hat{\G}_{12}\P {x_2}\cr
\Q {x_1,x_3}-\hat{\G}_{31}\P {x_1}-\hat{\G}_{13}\P {x_3}\cr
\Q {x_2,x_3}-\hat{\G}_{32}\P {x_2}-\hat{\G}_{23}\P {x_3}
\end{matrix}\right)
\end{equation}
where
$$
\begin{aligned}
&\hat{\G}_{21}=\frac{2(z-x)z^4}{(yz^3-1)q},\ \ \hskip 62 pt\hat{\G}_{12}=\frac{1-yz^3}{(z-x)q},\cr 
&\hat{\G}_{31}=\frac{2(z^4y-2yxz^3+2z-x)}{(yz^3-1)q},\ 
\ \hat{\G}_{13}=\frac{(2yz^3+2z^2+1)(1-yz^3)}{(z^4y-2yxz^3+2z-x)q},\cr 
&\hat{\G}_{32}=\frac{(z^4y-2yxz^3+2z-x)}{(x-z)zq},\ \ \ \ \hat{\G}_{23}=
\frac{(2z^2-4xz-3)(z-x)z^3}{(z^4y-2yxz^3+2z-x)q},
\end{aligned}
$$
and where
$$
q=2z^2+yz^3-2xz-1.
$$
These coefficients give rise to new Lam\'e potentials,
$$
h_1=\frac{yz^3-1}{q},\ \ h_2=\frac{(z-x)z}{q},\ \ h_3=\frac{2z-x-2yxz^3+z^4y}{q}.
$$  
The corresponding explicit solution of the 3-wave resonant interaction equations is
\begin{equation}\label{3WRIsol1}
\boxed{
\left(\begin{matrix}
0&A_{12}&A_{13}\cr
A_{21}&0&A_{23}\cr
A_{31}&A_{32}&0  
\end{matrix}\right)
=\frac{1}{q}\left(\begin{matrix} 0&-z & 2yz^3+2z^2+1\cr 
\ \ \   2z^3\ \  &\ \ \ 0 & \ \ \ z^2(2z^2-4xz-3)\cr 
\ \ \  2 &\ \ \  1 & 0\end{matrix}\right)
}
\end{equation}

We have shown that commencing with a certain diagonal action of $\mathbb{R}^4$ on $J^2(\mathbb{R},\mathbb{R})\times J^2(\mathbb{R},\mathbb{R})\times J^4(\mathbb{R},\mathbb{R})$ allows us to explicitly construct a solution of the {\sf m3WRI} equations and solution (\ref{3WRIsol0}) of the {\sf 3WRI} equations. Furthermore, because we are able to integrate differential system $H$ explicitly by methods of [\ref{AFV09}], we were able to construct, from solution (\ref{3WRIsol0}) of the {\sf 3WRI} system, the new solution (\ref{3WRIsol1}). 
Evidently, the formula (\ref{DMZsol}) permits the construction of a solution of the {\sf 3WRI} system that depends on 3 functions, each of one variable, $A(x),\, B(y),\, C(z)$. For different choices of these functions, these solutions all belong to the same gauge class and are, to some extent, labeled by the gauge invariants $h_{ij}$ of the {\sf DMZ} linear system (\ref{m3WRIsol}). In this case, four of the six gauge invariants vanish $h_{12}=h_{21}=h_{13}=h_{31}=0$ while
$$
h_{32}=\frac{2(y^2z^6-2yz^3+1)}{(x-yz^4+2xyz^3-2z)^2},\ \ h_{23}=\frac{6z^2(z-x)^2}{(x-yz^4+2xyz^3-2z)^2}.
$$ 
The fact that $h_{32}$ and $h_{23}$ are nonzero means that we can perform Laplace transformations in the $(3,2)$ and $(2,3)$-directions [\ref{KT96}] to obtain new {\sf DMZ} systems and consequently solutions of the {\sf 3WRI} and {\sf m3WRI} systems that belong to a gauge class distinct from that of solution (\ref{3WRIsol1}). 
\end{exmp}\hfill $\spadesuit$

\begin{rem}
Let $\mathbb{J}=J^{k_1}(\mathbb{R},\mathbb{R})\times J^{k_2}(\mathbb{R},\mathbb{R})\times J^{k_3}(\mathbb{R},\mathbb{R})$ and let
$$
\mu_D: \mathbb{R}^r\times\mathbb{J}\to\mathbb{J}
$$ 
be an admissible diagonal action of a Lie group $G$ as described in Theorems \ref{DMZdistributions1} and \ref{DMZdistributions2}, such that 
\begin{equation}\label{dimensionConstraint}
\dim\mathbb{J}-\dim G=10. 
\end{equation}
By choosing $k_i$ and $\dim G$ to be arbitrarily large, while maintaining the constraint (\ref{dimensionConstraint}) we can, by the means described above, obtain solutions of the {\sf 3WRI} system and the {\sf m3WRI} system, which depend on 3 smooth functions each of one variable and any number of their derivatives. An interesting question is whether or not these solutions are dense in the space of all solutions of the (say) {\sf 3WRI} system. 

Despite their importance, up to now very \textcolor{red}{\bf few} explicit {\sf DMZ} systems were known in the literature. Significant known {\sf DMZ} systems are firstly those arising from known triply orthogonal coordinate systems and the  ``1-periodic"  systems, derived by Kamran \& Tenenblat [\ref{KT98}]. Another example is given in [\ref{KT96}, Example 3]; however, the candidate system may contain a typographical error since it doesn't appear to be involutive.
\end{rem}

\section{Semi-Hamiltonian systems and $n$-hyperbolic manifolds}

A second application of our geometric construction of {\sf GDMZ} systems is related to strongly hyperbolic PDE systems of {\it hydrodynamic type} expressible in Riemann invariants in the form
\begin{equation}\label{hydroType}
u^i_t=v^i(u^1,u^2,\ldots,u^n)u^i_x,\ \ 1\leq i\leq n.
\end{equation}
The term {\it strongly hyperbolic} here means that the functions $v^i(u)$ in (\ref{hydroType}) are pair-wise distinct.

Here it is shown how to use the theory developed in previous sections for linear {\sf DMZ} systems of the form
\begin{equation}\label{linearKT}
\frac{\partial^2 u}{\partial x_i\partial x_j}-\Gamma^i_{ij}\frac{\partial u}{\partial x_i}-
\Gamma^i_{ij}\frac{\partial u}{\partial x_j}=0,\ 1\leq i<j\leq n
\end{equation}
to construct new examples of PDE systems (\ref{hydroType}) which are {\it semi-Hamiltonian}. Such systems play a role in a range of applications in mathematics and physics and have been the subject of interesting developments; for instance see [\ref{Ferapontov97}, \ref{FeraMoroSokolov09}]. The geometric study of such systems arose in the work of Dubrovin \& Novikov [\ref{DN}], D. Serre [\ref{Se}] and Tsarev [\ref{Tsa}] and we will briefly review some of their results.


There is an interesting subclass of strongly hyperbolic systems of hydrodynamic type which we now define.

\vskip 5 pt

\begin{defn}
A strongly hyperbolic system (\ref{hydroType}) is said to be {\it semi-Hamiltonian} or {\it rich in conservation laws} if
\begin{equation}\label{def_semiHam}
\left(\frac{v^i_{,j}}{v^j-v^i}\right)_{,k}=\left(\frac{v^i_{,k}}{v^k-v^i}\right)_{,j}
\end{equation}
\end{defn}



\begin{exmp} The so called {\it chromotography system}
$$
u^i_t=u^i\prod_{l=1}^Nu^lu^i_x,\ \ 1\leq i\leq N,
$$ 
is Hamiltonian when $N=2$ but only semi-Hamiltonian for $N>2$. 
\end{exmp} 
Various important results have been proven about semi-Hamiltonian systems. For instance, D. Serre [\ref{Se}] proved that Lax's classical result [\ref{Lax64}] on the blow-up in finite time for generic $2\times 2$ systems extends to the class of 
$n\times n$ semi-Hamiltonian systems. Another interesting result is due to Tsarev. To describe it we require a result of Darboux [\ref{Darboux10}]; see also  [\ref{Tenenblat98}].

The following paraphrases results of Tsarev [\ref{Tsa}].

\begin{thm}[Tsarev,{[\ref{Tsa}]}]\label{TsarevSemiham}
Let 
\begin{equation}\label{hydroType_thm}
u^i_t=v^i(u^1,u^2,\ldots,u^n)u^i_x,\ \ 1\leq i\leq n,
\end{equation}
be a strongly hyperbolic system of hydrodynamic type and consider the linear overdetermined PDE system
\begin{equation}\label{linearP_KT}
\frac{\partial^2 P}{\partial u_i\partial u_j}-\frac{v^i_{,j}}{v^j-v^i}\frac{\partial P}{\partial u_i}-
\frac{v^j_{,j}}{v^i-v^j}\frac{\partial P}{\partial u_j}=0,\ 1\leq i<j\leq n
\end{equation}
for a real-valued function $P(u^1,\ldots,u^n)$. Then (\ref{hydroType_thm}) is semi-Hamiltonian if and only if (\ref{linearP_KT}) is a {\sf DMZ} system.
\end{thm}

{\proof    Suppose (\ref{linearP_KT}) is a {\sf DMZ} system. The integrability conditions
\begin{equation}\label{KT_semihamCompatibility}
\begin{aligned}
&\frac{\pl\G^j_{ij}}{\pl u_k}-\G^k_{ik}\G^j_{kj}-\G^i_{ki}\G^j_{ij}+\G^j_{ij}\G^j_{kj}=0,\cr 
&\frac{\pl \G^j_{ij}}{\pl u_k}-\frac{\pl \G^j_{kj}}{\pl u_i}=0,
\end{aligned}(i,j,k)\in\ \text{\rm perm}\,(1,2,3).
\end{equation}
are therefore satisfied, where
\begin{equation}\label{Darboux}
\Gamma^i_{ij}(u)=\frac{v^i_{\,,j}}{v^j-v^i},\ 1\leq i\neq j\leq n.
\end{equation}
Equations $(\ref{linearP_KT})_2$ are easily seen to entail the condition (\ref{def_semiHam}). 

Conversely, if (\ref{hydroType_thm}) is semi-Hamiltonian, then a calculation verifies that (\ref{KT_semihamCompatibility}) will be satisfied given that $\G^i_{ij}(u)$  also satisfy (\ref{Darboux}).\hfill\qed 
}

\vskip 5 pt
In order to use the theory of this paper to construct new examples of semi-Hamiltonian hydrodynamic type systems, we are required from a given {\sf DMZ} system to solve the linear PDE system (\ref{Darboux}). Remarkably, the integrability condition arising from (\ref{Darboux}) is precisely $(\ref{KT_semihamCompatibility})_1$. This fact is attributed to Darboux [\ref{Darboux10}] in [\ref{Tenenblat98}], where a proof is given.

\begin{thm} [Darboux]\label{DarbouxTHM}
Let $\Gamma^i_{ij}(u^1,\ldots,u^n)$,\ $1\leq i\neq j\leq n$, be a collection of $n(n-1)$ smooth functions, symmetric in their lower indices and satisfying the linear {\sf DMZ} system integrability conditions with $C_{ij}=0$, equations (\ref{KT_0Compatibility}). That is,
\begin{equation}\label{KT_semihamCompatibilityDarboux}
\begin{aligned}
&\frac{\pl\G^j_{ij}}{\pl u_k}-\G^k_{ik}\G^j_{kj}-\G^i_{ki}\G^j_{ij}+\G^j_{ij}\G^j_{kj}=0, 
\end{aligned}(i,j,k)\in\ \text{\rm perm}\,(1,2,3).
\end{equation}
Then the linear first order PDE system 
\begin{equation}\label{Darboux2}
\frac{1}{w^j-w^i}\frac{\partial w^i}{\partial u^j}=\Gamma^i_{ij}(u),\ 1\leq i\neq j\leq n,
\end{equation}
has smooth solutions depending on $n$ arbitrary functions each of one variable.  
\end{thm}
\vskip 5 pt
Note that if we have by some means constructed a {\sf DMZ} system of the form
\begin{equation}\label{DMZconserved_density}
\frac{\pl^2 P}{\pl x_i\pl x_j}-\G^i_{ji}(u)\frac{\pl P}{\pl x_i}-\G^j_{ij}(u)\frac{\pl P}{\pl x_j}=0
\end{equation}
for a function $P(u^1,\ldots,u^n)$, then its coefficients will automatically satisfy (\ref{KT_semihamCompatibility}).

The significance of this remark arises from the following result of Tsarev [\ref{Tsa}]. 
\vskip 5 pt
\begin{thm}[\,Tsarev,\ {[\ref{Tsa}]}\,]\label{TsarevMain} 
Let 
\begin{equation}\label{hydroTypeTHM}
u^i_t=v^i(u^1,u^2,\ldots,u^n)u^i_x,\ \ 1\leq i\leq n
\end{equation}
be a semi-Hamiltonian system of hydroynamic type. Then
\begin{enumerate}
\item[1)] If $w^i(u)$ is a solution of 
\begin{equation}\label{semihamFromDMZ}
\frac{1}{w^j-w^i}\frac{\partial w^i}{\partial u^j}=\Gamma^i_{ij}(u)=\frac{1}{v^j-v^i}\frac{\partial v^i}{\partial u^j},\ 1\leq i\neq j\leq n,
\end{equation}
then 
$$
u^i_t=w^i(u)u^i_x
$$
is semi-Hamiltonian and defines a flow that commutes with (\ref{hydroTypeTHM}).
\item[2)] If  $w^i(u)$ is a solution of (\ref{semihamFromDMZ}) then  the vector field
$$
\sum_{i=1}^nw^i(u)u^i_x\P {u^i}
$$
is a first order generalised symmetry of (\ref{hydroTypeTHM}).
\end{enumerate}
\end{thm}

{\proof This follows Darboux's Theorem \ref{DarbouxTHM} and Theorem \ref{TsarevSemiham} }. \hfill\qed 

\vskip 5 pt

Additionally, suppose (\ref{hydroTypeTHM}) is semi-Hamiltonian with commuting flow components $w^i(u)$. Then Tsarev considers the $n$ algebraic equations
\begin{equation}\label{hodograph}
w^i(u)=v^i(u)t+x,\ \ 1\leq i\leq n,
\end{equation}
and proves
\begin{thm}[\,Tsarev,\ {[\ref{Tsa}]}\,] 
Suppose (\ref{hydroTypeTHM}) is semi-Hamiltonian with commuting flow components $w^i(u)$. A smooth solution $u^i(x,t)$ of the algebraic system (\ref{hodograph}) is a solution of (\ref{hydroTypeTHM}). Conversely, any solution $u^i(x,t)$ of (\ref{hydroTypeTHM}) can be locally represented as a solution of (\ref{hodograph}) in a neighbourhood of a point $(x_0,t_0)$ such that $u_x^i(x_0,t_0)\neq 0$ for every $i$.
\end{thm}

\begin{rem}
Note that functions $P$ on phase space that satisfy (\ref{linearP_KT}) are conserved densities of a semi-Hamiltonian system. Thus a semi-Hamiltonian system has infinitely many conserved densities.
\end{rem} 

Given the rich structure of semi-Hamiltonian systems descibed above it would seem to be important to find a means of constructing them. To date it appears that not many are known. However, the results of this paper will permit us to construct infinitely many examples. The idea is to use our construction of Darboux integrable $n$-hyperbolic manifolds, that is, Theorem \ref{DMZdistributions2}, to construct linear {\sf DMZ} systems from which we will be able to construct semi-Hamiltonian systems of hydrodynamic type and their commuting flows.
\vskip 10 pt
\subsection{Semi-Hamiltonian systems from 3-hyperbolic manifolds} The lowest dimension $n$ in which integrability conditions appear is the case $n=3$. As an application of the theory we will construct a simple example of a semi-Hamiltonian system of hydrodynamic type using an action of $\mathbb{R}^2$. It transpires that the associated Riemannian 3-metric is not flat and so our hydrodynamic system is not Hamiltonian. Nevertheless, we will exhibit an infinite family of commuting flows. 

In accordance with our general theory we can commence with the same manifold 
$J^2(\mathbb{R},\mathbb{R})\times J^2(\mathbb{R},\mathbb{R})\times J^2(\mathbb{R},\mathbb{R})$  as in Example \ref{exC22C2@A1} with the only difference that we quotient by an action of $G=\mathbb{R}^2$. As before set $\S_1=J^2(\mathbb{R},\mathbb{R})\times J^2(\mathbb{R},\mathbb{R})$ and $\S_2=J^2(\mathbb{R},\mathbb{R})$. Let 
$$
\boldsymbol{x}=(x,m,m_1,m_2,y,n,n_1,n_2)
$$
be coordinates on $\S_1$ and
$$
\boldsymbol{z}=(z,q,q_1,q_2)
$$
coordinates on $\S_2$. Let $\mu_1:G\times \S_1\to \S_1$ be the Lie group action
$$
\mu_1(g)(\boldsymbol{x})=(x,m-xt_1+t_2,m_1-t_1,m_2,y,n-yt_1+t_2,n_1-t_1,n_2).
$$
where $g=(t_1,t_2)$. By the standard procedure [\ref{FO}], one can take a right moving frame
$\rho:M_1\to G$ to be 
$$
\rho(\boldsymbol{x})=\left(\frac{m-n}{x-y},\frac{my-nx}{x-y}\right)
$$
The invariants of $\mu_1$ can be taken to be
$$
\left(x,\ m_1-\frac{m-n}{x-y},\ m_2,\ y,\ n_1-\frac{m-n}{x-y},\ n_2\right)
$$
These invariants together with the components of the moving frame form a local coordinate system on $\S_1$. Similarly, on $\S_2$ we define the action $\mu_2:G\times \S_2\to \S_2$
$$
\mu_2(g)(\boldsymbol{z})=(z,q+t_1z-t_1,q_1+t_1,q_2)
$$  
We can now form the map
$$
\boldsymbol{\phi^1}:\S_1\to\mathbb{R}^6\times G_{\boldsymbol{a}}
$$
by
$$
\begin{aligned}
\boldsymbol{\phi^1}(\boldsymbol{x},\boldsymbol{y})=&\left(x,\ m_1-\frac{m-n}{x-y},\ m_2,\ y,\ n_1-\frac{m-n}{x-y},\ n_2,\frac{m-n}{x-y},\frac{my-nx}{x-y}\right)\cr
                =&(x,\varpi_1,u_1,y,\varpi_2,u_2,a_1,a_2)
\end{aligned}
$$
and we get
$$
\begin{aligned}
\boldsymbol{\phi^1}_*\mathcal{C}\langle 2,2\rangle_{|_{\a=w}}=&\Big\{\P x+u_1\P {\varpi_1}-\frac{\varpi_1}{x-y}\left(\P {\varpi_1}+\P {\varpi_2}-\P {a_1}-y\P {a_2}\right),\cr
                                    &\hskip 20 pt\P y+u_2\P {\varpi_2}+\frac{\varpi_2}{x-y}(\P {\varpi_1}+\P {\varpi_2}-\P {a_1}-x\P {a_2}),\P {u_1},\P {u_2}\Big\}\cr
                                    =&\{\tilde{X},\tilde{Y},\P {u_1},\P {u_2}\}.
                                    \end{aligned}
$$
Similarly, the map
$$
\boldsymbol{\phi^2} : \S_2\to \mathbb{R}^2\times G_{\boldsymbol{b}}
$$
is
$$
\boldsymbol{\phi^2}(\boldsymbol{z})=(z,q_2,-q_1,q-zq_1)=(z,v_1,b_1,b_2)
$$
and we get
$$
\boldsymbol{\phi^2}_*\mathcal{C}\langle 2\rangle=\{\P z-v_1(\P {b_1}+z\P {b_2}),\P {v_1}\}=\{\tilde{Z},\P {v_1}\}.
$$
It can be verified that
$$
(\S_1\times_G\S_2,H)
$$
where
$$
H=\{X,\P {u_1}\}\oplus\{Y,\P {u_2}\}\oplus\{Z,\P {v_1}\}
$$
and 
$$
\begin{aligned}
&X=\P x+u_1\P {\varpi_1}-\frac{\varpi_1}{x-y}\left(\P {\varpi_1}+\P {\varpi_2}-\P {w_1}-y\P {w_2}\right),\cr
&Y=\P y+u_2\P {\varpi_2}+\frac{\varpi_2}{x-y}\left(\P {\varpi_1}+\P {\varpi_2}-\P {w_1}-x\P {w_2}\right),\cr
&Z=\P z-v_1(\P {w_1}+z\P {w_2}),
\end{aligned}
$$
is a Darboux integrable 3-hyperbolic manifold and we can use Theorem \ref{KTdistributions} to construct the associated {\sf DMZ} system.

It is easy to see that the invariants of $\mathcal{A}$ are 
$$
x,y,z,w_2-zw_1,
$$
and hence we compute as in Theorem \ref{KTdistributions} to get
\begin{equation}\label{PDEdefnEqs2_1}
p=w_2-zw_1,p_1=Xp=\frac{z-y}{x-y}\varpi_1,p_2=Yp=\frac{x-z}{x-y}\varpi_2,p_3=w_1
\end{equation}
and
\begin{equation}\label{PDEdefnEqs2_2}
p_{12}=\frac{z-x}{(x-y)^2}\varpi_1+\frac{z-y}{(x-y)^2}\varpi_2,p_{13}=\frac{\varpi_1}{x-y},p_{23}=-\frac{\varpi_2}{x-y}
\end{equation}
Solving equations (\ref{PDEdefnEqs2_1}) for $\varpi_1,\varpi_2,w_1,w_2$ in terms of $p,p_1,p_2,p_3$ we obtain the linear {\sf DMZ} system
\begin{equation}\label{KTeqsAbV2}
\begin{aligned}
&u_{xy}-\frac{x-z}{(x-y)(y-z)}u_x+\frac{y-z}{(x-y)(x-z)}u_y=0,\cr
&u_{xz}+\frac{1}{y-z}u_x=0,\ u_{yz}+\frac{1}{x-z}u_y=0.
\end{aligned}
\end{equation}
By Theorem 2.3, this system is involutive. We have
\begin{equation}\label{Christoffel1}
\begin{aligned}
&\G^1_{12}=\frac{x-z}{(x-y)(y-z)},\ \ \G^2_{12}=-\frac{y-z}{(x-y)(x-z)},\cr
&\G^1_{13}=-\frac{1}{y-z},\ \ \G^3_{13}=0,\cr
&\G^2_{23}=-\frac{1}{x-z},\ \ \G^3_{23}=0.
\end{aligned}
\end{equation}


A direct substitution of functions (\ref{Christoffel1}) into equations (\ref{DMZfullCompatibility}) verifies that system (\ref{KTeqsAbV2}) is involutive. It can further be shown that all the higher-dimensional Laplace invariants for (\ref{KTeqsAbV2}) vanish so an explicit solution for them can be constructed.

Now according to Darboux's theorem, Theorem \ref{DarbouxTHM}, the overdetermined system
$$
\begin{aligned}
&\frac{\partial w^1}{\partial u^2}=\frac{u^1-u^3}{(u^1-u^2)(u^2-u^3)}(w^2-w^1), \ \ \frac{\partial w^1}{\partial u^3}=-\frac{1}{u^2-u^3}(w^3-w^1),\cr
&\frac{\partial w^2}{\partial u^1}=-\frac{u^2-u^3}{(u^1-u^2)(u^1-u^3)}(w^1-w^2),\ \ \ \frac{\partial w^2}{\partial u^3}=-\frac{1}{u^1-u^3}(w^3-w^2),\cr
&\frac{\partial w^3}{\partial u^1}=0,\ \ \ \ \frac{\partial w^3}{\partial u^2}=0.
\end{aligned}
$$ 
is involutive and has solutions depending on three arbitrary functions each of one variable.  Indeed, it is not difficult to explicitly solve this linear system to get
\begin{equation}\label{ThreeComponentSemi}
\begin{aligned}
&w^1(u)=\frac{1}{u^3-u^2}\left((u^2-u^1)f_1'(u^1)+f_1(u^1)-f_2(u^2)+f_3(u^3)\right)\cr
&w^2(u)=\frac{1}{u^3-u^1}\left((u^2-u^1)f_2'(u^2)+f_1(u^1)-f_2(u^2)+f_3(u^3)\right)\cr
&w^3(u)=f_3'(u^3)
\end{aligned}
\end{equation}
In accordance with Theorem \ref{TsarevMain}, the system of hydrodynamic type
$$
u^i_t=w^i(u)u^i_x
$$
is semi-Hamiltonian for each solution $w(u)$ given by (\ref{ThreeComponentSemi}). Moreover, the collection of all such systems are commuting flows.

The linear system (\ref{KTeqsAbV2}) is but one of an infinite family of such systems, obtained by carrying out the same calculation but replacing the action $\mu^1$ by
$$
\begin{aligned}
&g\cdot (x,m,m_1,m_2,y,n,n_1,n_2)=\Big(x,\ m-h(x)t_1+t_2,\ m_1-h'(x)t_1,\ m_2-h''(x)t_1,\cr
&\hskip 150 pt y,\ n-k(y)t_1+t_2,\ n_1-k'(y)t_1,\ n_2-k''(y)t_1\Big)
\end{aligned}
$$
and replacing $\mu_2$ by
$$
g\cdot (z,p,p_1,p_2)=\Big(z,\ p+g(z)t_1-t_2,\ p_1+g'(z)t_1,\ p_2+g''(z)t_1\Big)
$$ 
where $h(x),k(y)$ and $g(z)$ are any $C^3$ functions. We thereby obtain the {\sf DMZ} system
\begin{equation}\label{2dAction_general}
\begin{aligned}
&u_{xy}-\frac{(g-h)k'}{(g-k)(h-k)}u_x+\frac{(g-k)h'}{(g-h)(h-k)}u_y=0,\cr
&u_{xz}-\frac{g'}{g-k}u_x=0,\cr
&u_{yz}-\frac{g'}{g-h}u_y=0.
\end{aligned}
\end{equation}
whose coefficients therefore satisfy the integrability conditions (\ref{DMZfullCompatibility}) and hence, by Theorem 7, determine semi-Hamiltonian systems of hydrodynamic type and their commuting flows. It can be checked that in this case, the Laplace invariants associated with system (\ref{2dAction_general}) turn out to be all zero and therefore, in some sense equation (\ref{2dAction_general}) is the {\it simplest} {\sf DMZ} system and the corresponding semi-Hamiltonian and {\sf 3WRI}/{\sf m3WRI} systems are the ``simplest" ones.

The examples given in this paper serve only to illustrate our geometric characterisation and construction of {\sf GDMZ} systems rather than possessing intrinsic significance of themselves. Their goal has been to provide a pointer to the possible new perspectives that the theory developed herein may provide. A much more detailed and complete analysis of semi-Hamiltonian systems of conservation laws as well as explicit exact solutions of 2+1-dimensional {\sf 3WRI} and {\sf m3WRI} systems will be given in forthcoming work.

\vskip 10 pt
\noindent{\bf Acknowledgements.}
I am indebted to my colleague Sergey Sergeev for numerous discussions on the theme of this paper and crucially, for pointing out the relationship between the three-wave resonant interaction systems (standard and modified) and the integrability conditions of the Darboux-Manakov-Zakharov linear problem.  On the basis of our discussions, several related works are in preparation on both continuous and discrete three-wave systems.

%
%
%
%

\end{document}